\documentclass[aap,MSNbibl,dvips]{arximspdf}
\usepackage{graphics}
\usepackage{mathrsfs}
\input xypic

\doi{10.1214/12-AAP915} 
\volume{24}
\issue{1}
\pubyear{2014}
\firstpage{76}
\lastpage{113}

\makeatletter
\def\mid{|}

\newcommand{\veps}{\varepsilon}

\newcommand{\diag}{\mathop{\operatorname{diag}}}

\newtheorem{theorem}{Theorem}[section]
\newtheorem{lemma}[theorem]{Lemma}
\newtheorem{corollary}[theorem]{Corollary}
\newtheorem{proposition}[theorem]{Proposition}
\newproclaim{definition}[theorem]{Definition}
\newproclaim{asmption}[theorem]{Assumption}

\newcommand{\cS}{\mathcal{S}}
\newcommand{\bLambda}{\bolds{\Lambda}}
\newcommand{\bGamma}{\bolds{\Gamma}}
\newcommand{\bpi}{\bolds{\pi}}
\newcommand{\brho}{\bolds{\rho}}
\newcommand{\bmu}{\bolds{\mu}}
\newcommand{\bM}{\mathbf{M}}
\newcommand{\bnu}{\bolds{\nu}}
\newcommand{\bkappa}{\bolds{\kappa}}
\newcommand{\btheta}{\bolds{\theta}}
\newcommand{\bseta}{\bolds{\eta}}
\newcommand{\bA}{\mathbf{A}}
\newcommand{\bC}{\mathbf{C}}

\newcommand{\Rpp}{\mathbb{R}_{+}}
\newcommand{\R}{\mathbb{R}}
\newcommand{\Zp}{\mathbb{Z}_{+}}
\newcommand{\cJ}{\mathcal{J}}
\newcommand{\cI}{\mathcal{I}}

\newcommand{\E}{\mathbb{E}}
\newcommand{\PP}{\mathbb{P}}

\newcommand{\bQ}{\mathbf{Q}}
\newcommand{\bN}{\mathbf{N}}
\newcommand{\bW}{\mathbf{W}}

\newcommand{\by}{\mathbf{y}}
\newcommand{\bbm}{\mathbf{m}}
\newcommand{\be}{\mathbf{e}}
\newcommand{\bv}{\mathbf{v}}

\newcommand{\bone}{\mathbf{1}}

\newcommand{\eqref}[1]{(\ref{#1})}
\makeatother

\begin{document}
\begin{frontmatter}

\title{Qualitative properties of $\alpha$-fair policies in
bandwidth-sharing networks}
\runtitle{Qualitative properties of $\alpha$-fair policies}

\begin{aug}
\author[a]{\fnms{D.} \snm{Shah}\ead[label=e1]{devavrat@mit.edu}},
\author[a]{\fnms{J.~N.} \snm{Tsitsiklis}\ead[label=e2]{jnt@mit.edu}}
\and
\author[b]{\fnms{Y.} \snm{Zhong}\corref{}\thanksref{ds.ack}\ead[label=e3]{zhyu4118@berkeley.edu}}
\runauthor{D. Shah, J.~N. Tsitsiklis and Y. Zhong}
\affiliation{Massachusetts Institute of Technology, Massachusetts
Institute of Technology and University of California}
\address[a]{D. Shah\\
J.~N. Tsitsiklis\\
Department of EECS\\
Massachusetts Institute of Technology\\
Cambridge, Massachusetts 02139\\
USA\\
\printead{e1}\\
\phantom{E-mail:\ }\printead*{e2}}

\address[b]{Y. Zhong\\
Department of EECS\\
University of California\\
Berkeley, California 94720\\
USA\\
\printead{e3}}
\thankstext{ds.ack}{Supported by NSF project CCF-0728554,
and was performed while the authors were affiliated with the Laboratory
for Information and Decision
Systems as well as the Operations Research
Center at MIT. The third author is now affiliated with the Algorithms,
Machines and People's
Laboratory at the University of California, Berkeley.} 
\end{aug}

\received{\smonth{4} \syear{2011}}
\revised{\smonth{10} \syear{2012}}

%
\begin{abstract}
We consider a flow-level model of a network operating under an $\alpha$-fair
bandwidth
sharing policy (with $\alpha> 0$) proposed by Roberts and Massouli\'{e}
[\textit{Telecomunication Systems} \textbf{15} (2000) 185--201].
This is a probabilistic model that captures the long-term aspects of
bandwidth sharing between users or flows in a communication network.

We study the transient properties as well as the
steady-state distribution of the model.
In particular, for $\alpha\geq1$, we obtain bounds on the maximum
number of flows
in the network over a given time horizon, by means of
a maximal inequality derived from the standard Lyapunov drift
condition. As a
corollary, we establish the full state space collapse property for all
$\alpha\geq1$.

For the steady-state distribution, we obtain explicit exponential tail
bounds on the number
of flows, for any $\alpha> 0$, by relying on a norm-like Lyapunov function.
As a corollary, we establish the validity
of the diffusion approximation developed by Kang et al.
[\textit{Ann. Appl. Probab.} \textbf{19} (2009) 1719--1780], in steady
state, for the
case where $\alpha= 1$ {and under a local traffic condition.}
\end{abstract}

%
\begin{keyword}[class=AMS]
\kwd[Primary ]{60K20}
\kwd{68M12}
\kwd[; secondary ]{68M20}
\end{keyword}

\begin{keyword}
\kwd{Bandwidth-sharing policy}
\kwd{state space collapse}
\kwd{$\alpha$-fair}
\kwd{heavy traffic}
\end{keyword}

\end{frontmatter}

\section{Introduction}\label{sec:model}
We consider a flow-level model of a network that operates under an
$\alpha$-fair bandwidth-sharing policy, and establish a variety of new
results on the resulting performance. These results include tail bounds
on the size of a maximal excursion during a finite time interval,
finiteness of expected queue sizes, exponential tail bounds under the
steady-state distribution, and the validity of the heavy-traffic
diffusion approximation in steady state.
We note that our results are to a great extent parallel and
complementary to our work on packet-switched networks, which was
reported in \cite{STZ10}.

In the remainder of this section, we put our results in perspective by
comparing them with earlier work, and conclude with some more details
on the nature of our contributions.

\subsection{Background}\label{sec:back}
The flow-level network model that we consider
was introduced by Roberts and Massouli\'e \cite{RM00}
to study the dynamic behavior of Internet flows. It builds on a static version
of the model that was proposed earlier by Kelly, Maulloo, and Tan \cite{KMT98},
and subsequently
generalized by Mo and Walrand \cite{MW00} who introduced 
a class of ``fair'' bandwidth-sharing policies parameterized by
$\alpha> 0$.

The most basic question regarding flow-level models concerns
necessary and sufficient conditions for stability, that is, for the
existence of a steady-state distribution for
the associated Markov process. This question was answered by
Bonald and Massouli\'e \cite{BM01} for the case of $\alpha$-fair
policies with $\alpha> 0$,
and
by de Veciana et al. \cite{DeV01} for the case of max-min fair policies
($\alpha\to\infty$)
and proportionally fair ($\alpha=1$) policies. In all cases, the
stability conditions turned out to be the natural deterministic
conditions based on mean arrival and service rates.

Given these stability results, the natural next question is whether the
steady-state expectation of the number of flows in the system is finite
and if so, to identify some nontrivial upper bounds. When $\alpha\geq1$,
the finiteness question can be answered in the affirmative, and
explicit bounds can be obtained by exploiting the same Lyapunov drift
inequality that had been used in earlier work to establish stability.
However, this approach does not seem to apply to the case where $\alpha
\in(0,1)$, which remains an open problem; this is one of the problems
that we settle in this paper.

A more refined analysis of the number of flows present in the system
concerns exponentially decaying bounds {on} the tail of its
steady-state distribution. We provide results of this form, together
with explicit bounds for the associated exponent. While a result of
this type was not previously available, we
take note of
related recent results by Stolyar \cite{S-LDP08} and Venkataramanan
and Lin 
\cite{VL09} who provide a precise asymptotic characterization of the %
exponent of the tail probability, in steady state, for the case
of switched networks (as opposed to flow-level network models). [To be
precise, their results concern the $(1+\alpha)$ norm of the vector of
flow counts
under maximum weight or pressure policies parameterized by $\alpha
>0$.] We believe that
their methods extend to the model considered here, without much
difficulty. However, their approach leads to a variational
characterization that
appears to be difficult to evaluate (or even bound) explicitly.
We also take note of work by Subramanian \cite{vijay10}, 
who establishes a large deviations principle for a class of
switched network models under maximum weight or pressure
policies with $\alpha= 1$.

The analysis of the steady-state distribution for underloaded networks
provides only partial insights about the transient behavior of the associated
Markov process. As an alternative, the heavy-traffic (or diffusion)
scaling of the network can lead to parsimonious approximations for the
transient behavior. A general
two-stage program for developing such 
diffusion approximations has been put forth by Bramson \cite{bramson}
and Williams \cite{williams}, and has been carried out in detail for
certain particular classes of
queuing network models. To carry out this program, one needs to: (i)~{provide} a detailed analysis of a related fluid model
when the network is critically loaded; and, (ii)~identify a unique
distributional limit of the associated diffusion-scaled processes by
studying a related
Skorohod problem. The first stage of the program was carried out by
Kelly and Williams \cite{KW04} 
who identified the invariant manifold of the associated
critically loaded fluid model. This further led to the proof by Kang et
al. \cite{KKLW09} of a multiplicative
state space collapse property, similar to results by Bramson \cite
{bramson}. We note that the above summarized results
hold under $\alpha$-fair policies with an arbitrary $\alpha> 0$.
The second stage of the program has been carried out for the proportionally
fair policy ($\alpha= 1$) by Kang et al. \cite{KKLW09}, 
under a technical \textit{local traffic} condition, and
more recently, by Ye and Yao \cite{YY10}, 
under a somewhat less restrictive 
technical condition.
{We note, however, that when $\alpha\neq1$, a diffusion approximation
has not been established. In this case, it is of interest to see at
least whether
properties that are stronger than multiplicative state space collapse
can be derived,
something that is accomplished in the present paper.}

The above outlined diffusion approximation results involve rigorous
statements on the finite-time behavior of the original process.
Kang et al. \cite{KKLW09}
further established that for the particular setting that they consider,
the resulting diffusion approximation
has an elegant product-form steady-state distribution; this result
gives rise to an
intuitively appealing interpretation of the relation between the
congestion control protocol utilized by the flows (the end-users)
and the queues formed inside the network.
It is natural to expect that this product-form steady-state
distribution is the limit
of the steady-state distributions in the original
model under the diffusion scaling. Results of this type are known for
certain queueing systems such as generalized Jackson
networks; see the work by Gamarnik and Zeevi \cite{GZ06}. 
On the other hand, the validity of
such a steady-state diffusion approximation
was not known for the model considered in \cite{KKLW09}; it will be
established in the present paper.

\subsection{Our contributions}
In this paper, we advance
the performance analysis of flow-level models
of networks operating under an $\alpha$-fair policy, in both the steady-state
and the transient regimes.


For the transient regime, we obtain a probabilistic bound on the
maximal (over a finite time horizon) number of flows,
when operating under an $\alpha$-fair policy with $\alpha\geq1$.
This result is obtained by combining a Lyapunov drift inequality with
a natural extension of Doob's maximal inequality
for nonnegative supermartingales. Our probabilistic bound, together
with prior results on
multiplicative state space collapse, leads immediately to a stronger
property, namely,
full state space
collapse, for the case where $\alpha\geq1$.

For the steady-state regime, we obtain
nonasymptotic and explicit bounds on the tail of the distribution of
the number of flows, for
any $\alpha> 0$. In the process, we establish that, for any $\alpha
>0$, all moments
of the steady-state number of flows are finite.
These results
are proved by working with a \textit{normed} version of the Lyapunov function
that was used in prior work. Specifically, we establish that
this normed version is also a Lyapunov function for the system (i.e.,
it satisfies a drift inequality).
It also happens to be a Lipschitz continuous function and this helps
crucially in establishing exponential tail bounds, using results of
Hajek \cite{Hajek-delay} and
Bertsimas et al. \cite{BGT01}.

The exponent in the exponential tail bound that we establish
for the distribution of the number of flows is proportional to a
suitably defined distance (``gap'') from critical loading; this gap is
of the same type as the familiar $1-\rho$ term, where $\rho$ is the
usual load factor in a queueing system. This particular dependence on
the load leads to the tightness
of the steady-state distributions of the model under diffusion
scaling. It leads to one of our main results, namely, the validity of
the diffusion approximation, in steady state, when $\alpha=1$ and a
local traffic condition holds.


\subsection{Organization}
The rest of the paper is organized as follows. In Section \ref{sec2}, we define
the notation
and some of the terminology that we will employ. We also
describe the flow-level network model, as well as the weighted
$\alpha$-fair bandwidth-sharing policies. In Section \ref{sec3}, we provide formal
statements of our main results.
The transient analysis is presented in Section \ref{sec4}. We start with a general
lemma, and specialize it to obtain a maximal inequality under $\alpha
$-fair policies, when
$\alpha\geq1$. We then apply the latter inequality to prove full
state space
collapse when $\alpha\geq1$.

We then proceed to the steady-state analysis.
In Section \ref{sec:drift}, we establish a drift inequality for a
suitable Lyapunov function, which is central to our proof of
exponential upper
bounds on tail probabilities 
We prove the exponential upper bound on tail probabilities in
Section \ref{sec:exp}. The validity of the heavy-traffic steady-state approximation
is established
in Section \ref{sec:inlimit}. 
We conclude the paper with a brief discussion in Section \ref{sec:concl}.

\section{Model and notation}\label{sec2}
\subsection{Notation}
We introduce here the notation
that will be employed throughout the paper.
We denote the real vector space of dimension\vadjust{\goodbreak} $M$ by $\mathbb{R}^M$,
the set of
nonnegative $M$-tuples by $\mathbb{R}_{+}^M$
and the set of positive $M$-tuples by $\mathbb{R}_{p}^M$. We write
$\mathbb{R}$ for $\mathbb{R}^1$,
$\mathbb{R}_{+}$ for $\mathbb{R}_+^1$ and $\mathbb{R}_{p}$ for
$\mathbb{R}_{p}^1$.
We let $\mathbb{Z}$ be the set of integers,
$\mathbb{Z}_+$ the set of nonnegative integers and $\mathbb{N}$ the
set of positive
integers. Throughout the paper, we reserve bold letters for vectors and
plain letters for scalars.

For any vector $\mathbf{x}\in\mathbb{R}^M$, and any $\alpha>0$, we define
\[
\| \mathbf{x}\|_{\alpha}= \Biggl( \sum_{i=1}^M
|x_i|^{\alpha} \Biggr)^{1/\alpha},
\]
and we define $\|\mathbf{x}\|_{\infty} = \max_{i \in\{1,\ldots,M\}
} |x_i|$.
For any two vectors $\mathbf{x}=(x_i)_{i=1}^M$
and $\mathbf{y}=(y_i)_{i=1}^M$ of the same {dimensions}, we let $
\langle{\mathbf{x},\by} \rangle=
\sum_{i=1}^M x_i y_i$ be the inner product of $\mathbf{x}$ and
$\mathbf{y}$.
We let $\mathbf{e}_i$ 
be the $i$th unit vector in
$\mathbb{R}^M$, and $\bone$ the vector of all ones. For a set
$\mathcal{S}$, we denote its
cardinality by $|\mathcal{S}|$, and its indicator function by $
\mathbb{I}_{\mathcal{S}}$. For a matrix $\bA$, we
let $\bA^T$ denote its transpose.

\subsection{Flow-level network model}\label{ssec:model}
\subsubsection*{The model}
We adopt the model and notation in \cite{KW04}. As explained in detail
in \cite{KW04}, this model faithfully captures the long-term (or macro
level) behavior of congestion control in the current Internet.

Let time be continuous and indexed
by $t \in\Rpp$. Consider a network with a finite set $\cJ$ of resources
and a set $\cI$ of routes, where a route is identified with a
nonempty subset
of the resource set $\cJ$. Let $\bA$ be the $|\cJ|\times|\cI|$
matrix with $A_{ji}=1$ if
resource $j$ is used by route $i$, and $A_{ji}=0$ otherwise. Assume
that $\bA$ has rank $|\cJ|$. Let
$\bC= (C_j)_{j\in\cJ}$ be a \textit{capacity} vector, where we assume that
each entry $C_j$ is a given positive constant. Let the number
of flows on route $i$ at time $t$ be denoted by $N_i(t)$, and define
the flow vector at time $t$ by $\bN(t)=(N_i(t))_{i\in\cI}$.
For each route $i$, new flows arrive
as an independent Poisson process of rate $\nu_i$. Each arriving flow
brings an amount of
work (data that it wishes to transfer) which is an exponentially
distributed random variable with
mean $1/\mu_i$, independent of everything else. Each flow
gets service from the network according to a bandwidth-sharing policy. Once
a flow is served, it departs the network.\vspace*{-2pt}

\subsubsection*{The $\alpha$-fair bandwidth-sharing policy} A bandwidth
sharing policy has to allocate rates to flows so that capacity
constraints are satisfied at each time instance. Here we discuss the
popular $\alpha$-fair bandwidth-sharing policy, where $\alpha> 0$. At
any time, the bandwidth allocation depends on the current number of
flows $\mathbf {n}=(n_i)_{i\in\cI}$. Let $\Lambda_i$ be the
\textit{total} bandwidth allocated to route $i$ under the $\alpha$-fair
policy: each flow of type $i$ gets rate $\Lambda_i/n_i$ if $n_i>0$, and
$\Lambda_i = 0$ if $n_i = 0$. Under an $\alpha$-fair policy, the
bandwidth vector $\bLambda(\mathbf{n}) = (\Lambda_i(\mathbf
{n}))_{i\in\cI}$ is determined as follows.

If $\mathbf{n}=\mathbf{0}$, then $\bLambda= \mathbf{0}$. If
$\mathbf{n}\neq\mathbf{0}$, then let
$\cI_{+}(\mathbf{n})=\{i\in\cI\dvtx  n_i>0\}$. For $i \notin\cI
_{+}(\mathbf{n})$, set $\Lambda_i(\mathbf{n})=0$.
Let $\bLambda_{+}(\mathbf{n})=(\Lambda_i(\mathbf{n}))_{i\in\cI
_{+}(\mathbf{n})}$.\vadjust{\goodbreak} Then{} $\bLambda_{+}(\mathbf{n})$ is
the \textit{unique} maximizer in the optimization problem
%
\begin{eqnarray}
&&\mathsf{maximize}\ G_n(\bLambda_{+})\ \mathsf{over}\
\bLambda \in\Rpp^{|\cI|} \label{utility}
\\
&&\mathsf{subject\ to} \qquad \sum_{i\in\cI_{+}(n)}A_{ji}
\Lambda_i\leq C_j, \forall j \in\cJ, \label{constraint1}
\end{eqnarray}
where
\[
G_{\mathbf{n}}(\bLambda_{+}) = \cases{ %
\displaystyle\sum_{i\in\cI_{+}(\mathbf{n})}\kappa_i
n_i^{\alpha
}\frac{\Lambda_i^{1-\alpha}}{1-\alpha}, & \quad$\mbox{if } \alpha\in (0,
\infty)\setminus\{1\},$
\vspace*{2pt}\cr
\displaystyle\sum_{i\in\cI_{+}(\mathbf{n})}\kappa_i n_i
\log \Lambda_i, &\quad $\mbox{if } \alpha= 1.$}
\]
Here, for each $i \in\cI$, $\kappa_i$ is a positive weight assigned
to route $i$.

Some crucial properties of $\bLambda(\mathbf{n})$ are as follows (see
Appendix A of \cite{KW04}):
\begin{longlist}[(iii)]
\item[(i)] $\Lambda_i(\mathbf{n})>0$ for every $i\in\cI
_{+}(\mathbf{n})$;
\item[(ii)] $\bLambda(r\mathbf{n})=\bLambda(\mathbf{n})$ for $r>0$;
\item[(iii)] For every $\mathbf{n}$ and every $i\in\cI_{+}(\mathbf{n})$,
the function $\Lambda_i(\cdot)$ is continuous at~$\mathbf{n}$.
\end{longlist}

\subsubsection*{Flow dynamics} The flow dynamics are described by the
evolution of the flow vector $\bN(t) = (N_i(t))_{i\in\cI}$, a Markov
process with {infinitesimal} transition rate matrix $\mathbf{q}$ given
by
%
\begin{equation}
\label{mx:tsn} q(\mathbf{n},\mathbf{n}+\bbm) = \cases{ %
\nu_i, & \quad$\mbox{if }\bbm=\be_i,$
\vspace*{2pt}\cr
\mu_i\Lambda_i(\mathbf{n}), & \quad$\mbox{if }\bbm= -
\be_i \mbox{ and } n_i\geq1,$
\vspace*{2pt}\cr
0, &\quad $\mbox{otherwise,}$}
\end{equation}
where for each $i$, $\nu_i>0$ and $\mu_i>0$ are the arrival and
service rates defined earlier,
and $\be_i$ is the $i$th unit vector.

\subsubsection*{Capacity region} Flows of type $i$ bring to the system
an average of $\rho_i = \nu_i/\mu_i$ units of work per unit time.
Therefore, in order for the Markov process $\bN(\cdot)$ to be positive
recurrent, it is necessary that
%
\begin{equation}
\label{eq:capfl}
\bA\brho <  \bC\qquad \mbox{componentwise}.
\end{equation}
We note that under the $\alpha$-fair bandwidth-sharing policy,
condition \eqref{eq:capfl} is also sufficient for
positive recurrence of the process $\bN(\cdot)$ \cite{BM01,DeV01,KW04}.



\subsection{A note on our use of constants}
\label{ssec:rmk}
Our results and proofs involve various constants; some are absolute
constants, some depend only on the structure of the network and some
depend (smoothly) on the traffic parameters (the arrival and service
rates). It is convenient to distinguish between the different types of
constants, and we define here the terminology that we will be
using.\vadjust{\goodbreak}

The term \textit{absolute constant} will be used to refer to a quantity
that does not depend on any of the model parameters. The term {{\em
network-dependent constant}} will be used to refer to quantities that
are completely determined by the structure of the underlying network
and policy, namely,
the incidence matrix $\bA$, the capacity vector $\bC$, the weight
vector $\bkappa$
and the policy parameter $\alpha$.

Our analysis also involves certain quantities that depend on the
traffic parameters, namely, the arrival and service parameters $\bmu$
and $\bnu$. These quantities are often given by complicated
expressions that would be inconvenient to carry through the various
arguments. It turns out that the only property of such quantities that
is relevant to our purposes is the fact they change continuously as
$\bmu$ and $\bnu$ vary over the open positive orthant.
(This still allows these quantities to be undefined or discontinuous on
the boundary of the positive orthant.) We abstract this property by
introducing, in the definition that follows, the concept of a
{\emph{{(positive)} load-dependent constant}}.

%
\begin{definition}\label{df:ldc}
Consider a family of bandwidth-sharing networks with common parameters
$(\bA, \bC, \bkappa,\alpha)$, but varying traffic parameters
$(\bmu,\bnu)$.\break
A~quantity~$K$ will be called a \textit{{(positive)}
load-dependent constant} if for networks in that family it is
determined by a relation of the form
$K = f(\bmu, \bnu)$, where $f\dvtx \R_p^{|\cI|}\times\R_p^{|\cI|}
\to \R_p$ is a \textit{continuous} function on the open positive
orthant $\R_p^{|\cI|}\times\R_p^{|\cI|}$.
\end{definition}

A key property of a load-dependent constant, which will be used in some
of the subsequent proofs, is that it is
by definition positive and furthermore (because of continuity),
bounded above and below by positive network-dependent constants if we
restrict $\bmu$ and $\bnu$ to a compact subset of the open positive orthant.
A natural example of a load-dependent constant is the load factor $\rho
_i=\nu_i/\mu_i$. (Note that this quantity diverges as $\mu_i\to0$.)

{We also define}
the \textit{gap} of a underloaded bandwidth-sharing network.

%
\begin{definition}
Consider a family of bandwidth-sharing networks with
common parameters $(\bA, \bC, \bkappa,\alpha)$ and with varying
traffic parameters $(\bmu,\bnu)$ that satisfy
$\bA\brho< \bC$. The \textit{gap} of a network with traffic
parameters $(\bmu,\bnu)$ in the family,
denoted by $\varepsilon(\brho)$, is defined by
\[
\varepsilon(\brho) \triangleq\sup \bigl\{\tilde{\varepsilon} > 0 \dvtx (1 +
\tilde{ \varepsilon}) \bA\brho\leq\bC \bigr\}.
\]
\end{definition}
%
We sometimes write $\varepsilon$ for $\varepsilon(\brho)$ when there
is no ambiguity.
Note also that $\varepsilon(\brho)$ plays the same role as the term
$1-\rho$ in a queueing system with load~$\rho$. 

\subsection{Uniformization} \label{se:unif}
Uniformization is a well-known device which allows us to study a
continuous-time Markov process by considering an associated
discrete-time Markov chain with the same stationary distribution. We
provide here some details and the notation that we will be
using.\vadjust{\goodbreak}

Recall that the Markov process $\bN(\cdot)$ of interest has dynamics
given by \eqref{mx:tsn}.
Let $\Xi(\mathbf{n}) = \sum_{\tilde{\mathbf{n}}} q(\mathbf
{n},\tilde{\mathbf{n}})$ be the aggregate transition rate at state
$\mathbf{n}$.
The \textit{embedded jump chain} of $\bN(\cdot)$ is a discrete-time
Markov chain with
the same state space $\Zp^{|\cI|}$, and with transition probability
matrix $\mathbf{P}$ given by
\[
P(\mathbf{n}, \tilde{\mathbf{n}}) = \frac{q(\mathbf{n},\tilde
{\mathbf{n}})}{\Xi(\mathbf{n})}.
\]
The so-called \textit{uniformized Markov chain} is an alternative, more
convenient, discrete-time Markov chain, denoted $ (\tilde{\bN
}(\tau) )_{\tau\in\Zp}$, to be defined shortly.

We first introduce some more notation.
Consider the aggregate transition rates $\Xi(\mathbf{n}) = \sum_{\tilde{\mathbf{n}}} q(\mathbf{n},\tilde{\mathbf{n}})$. Since
every route uses at least one resource, we have $\Lambda_i(\mathbf
{n}) \leq\max_{j\in\cJ} C_j$,
for all $i\in\cI$. Then, by \eqref{mx:tsn},
we have
\[
\Xi(\mathbf{n}) = \sum_{\tilde{\mathbf{n}}} q(\mathbf{n},\tilde {
\mathbf{n}}) \leq\sum_{i\in\cI} \bigl(\nu_i +
\mu_i \Lambda_i(\mathbf {n}) \bigr) \leq\sum
_{i\in\cI} \Bigl(\nu_i + \mu_i \max
_{j\in
\cJ} C_j \Bigr).
\]
We define $\Xi\triangleq\sum_{i\in\cI}  (\nu_i + \mu_i \max_{j\in\cJ} C_j )$, and modify the rates of self-transitions
(which were zero in the original model) to
%
\begin{equation}
\label{eq:tsn_update} q(\mathbf{n},\mathbf{n}):= \Xi- \Xi(
\mathbf{n}).
\end{equation}
Note that $\Xi$ is a positive load-dependent constant.
We define a transition probability matrix $\tilde{\mathbf{P}}$ by
\[
\tilde{P}(\mathbf{n},\tilde{\mathbf{n}}) \triangleq\frac
{q(\mathbf{n},\tilde{\mathbf{n}})}{\Xi}.
\]

\begin{definition}\label{df:umc}
The uniformized Markov chain $ (\tilde{\bN}(\tau) )_{\tau
\in\Zp}$
associated with the Markov process $\bN(\cdot)$ is a discrete-time
Markov chain
with the same state space $\Zp^{|\cI|}$, and with transition matrix
$\tilde{\mathbf{P}}$ defined as above.
\end{definition}


As remarked earlier, the Markov process $\bN(\cdot)$ that describes a
bandwidth-sharing network operating
under an $\alpha$-fair policy is positive recurrent, as long as the
system is underloaded,
that is, if $\bA\brho< \bC$. It is not hard to verify that $\bN
(\cdot)$ is also irreducible. Therefore,
the Markov process $\bN(\cdot)$ has a unique stationary distribution.
The chain $\tilde{\bN}(\cdot)$ is
also positive recurrent and irreducible because $\bN(\cdot)$ is, and
by suitably increasing $\Xi$
if necessary, it can be made aperiodic. Thus $\tilde{\bN}(\cdot)$
has a unique stationary
distribution as well.
A crucial property of the uniformized chain $\tilde{\bN}(\cdot)$ is
that this unique stationary distribution
is the same as that of the original Markov process $\bN(\cdot)$; see,
for example, \cite{dsp}.

\subsection{A mean value theorem}
We will be making extensive use of
a second-order mean value
theorem \cite{nlp}, which we state below for easy reference.

\begin{proposition}\label{prp:mvt}
Let $g \dvtx \mathbb{R}^M \rightarrow\mathbb{R}$ be twice continuously
differentiable over an open sphere $S$ centered at a vector $\mathbf{x}$.
Then, for any $\mathbf{y}$\vadjust{\goodbreak} such that $\mathbf{x}+\mathbf{y} \in S$, there
exists $\theta\in[0,1]$ such that
%
\begin{equation}
\label{eq:taylor} g(\mathbf{x}+\mathbf{y})=g(\mathbf{x})+\mathbf{y}^T
\nabla g(\mathbf{x})+ \tfrac{1}{2}\mathbf{y}^T H(\mathbf{x}+
\theta\mathbf{y}) \mathbf{y},
\end{equation}
where $\nabla g(\mathbf{x}) \triangleq [\frac{\partial
g(\mathbf{x})}{\partial x_i} ]_{i=1}^M \in\mathbb{R}^M$
is the gradient of $g$ at $\mathbf{x}$, and
\[
H(\mathbf{x}) \triangleq
\biggl[\frac{\partial^2 g(\mathbf{x})}{\partial x_i\, \partial
x_j} \biggr]_{i,j=1}^M
\in\mathbb{R}^{M\times M}
\]
is the Hessian of the function $g$ at $\mathbf{x}$.
\end{proposition}

\section{Summary of results}\label{sec3}
In this section, we summarize our main results 
for both the transient and the
steady-state regime. The proofs are given in subsequent sections.

\subsection{Transient regime}
Here we provide a simple inequality on the maximal excursion of the
number of flows over a
finite time interval, under an $\alpha$-fair policy with $\alpha\geq1$.
%
\begin{theorem}\label{thm:max}
Consider a bandwidth-sharing network operating under an $\alpha$-fair policy
with $\alpha\geq1$, and assume that $\bA\brho< \bC$.
Suppose that $\bN(0)=\mathbf{0}$. Let $N^*(T) = \sup_{t\in
[0,T], i\in\cI} N_i(t)$, and
let $\veps$ be the gap. 
Then, for any $b>0$,
%
\begin{equation}
\mathbb{P} \bigl(N^*(T)\geq b \bigr)\leq\frac{KT}{\veps^{\alpha
-1}b^{\alpha+1}}
\end{equation}
for some positive load-dependent constant $K$. 
\end{theorem}
As an important application, in Section \ref{ssec:ssc},
we will use Theorem \ref{thm:max} to prove a full state space collapse
result, 
when $\alpha\geq1$. {(As discussed in the \hyperref[sec:model]{Introduction},
this property is stronger than multiplicative state space collapse.)}
The precise statement can be found in Theorem \ref{thm:fssc}.

\subsection{Stationary regime}
As noted earlier, 
the Markov process $\bN(\cdot)$ has a unique stationary distribution,
which we will denote by $\bpi$.
We use $\E_{\bpi}$ and $\mathbb{P}_{\bpi}$
to denote expectations and probabilities under $\bpi$.


\subsubsection*{Exponential bound on tail probabilities}
For an $\alpha$-fair policy, and for any $\alpha\in(0,\infty)$, we
obtain an explicit exponential upper bound on the tail probabilities
for the number of flows, in steady state. This will be used to
establish an ``interchange of limits'' result in Section
\ref{sec:inlimit}. See Theorem \ref{thm:inlimit} for more details.
%
\begin{theorem}\label{thm:sw2}
Consider a bandwidth-sharing network operating under an $\alpha$-fair policy
with $\alpha> 0$,
and assume that $\bA\brho< \bC$. Let $\veps$ be the gap. 
There exist positive constants
$B$, $K$ and $\xi$ such that for all $\ell\in\Zp$,
%
\begin{equation}
\label{eq:exbdd} \mathbb{P}_{\bpi}\bigl (\|\bN\|_{\infty}\geq B + 2 \xi
\ell \bigr) \leq \biggl(\frac{\xi}{\xi+\veps K} \biggr)^{\ell+1}.\vadjust{\goodbreak}
\end{equation}
Here $\xi$ and 
$K$ 
are load-dependent constants,
and $B$ takes the form $K'/\veps$ when $\alpha\geq1$, and $K'/\min\{
\veps^{1/\alpha}, \veps\}$
when $\alpha\in(0,1)$, with $K'$ being a positive load-dependent constant.
In particular, all moments of $\|\bN\|_{\infty}$ are finite under the
stationary distribution $\bpi$,
that is, $\E_{\bpi}[\|\bN\|_{\infty}^k] < \infty$ for every $k \in
\mathbb{N}$.
\end{theorem}
Here we note that Theorem \ref{thm:sw2} implies the following.
The system load $L(\brho)$, defined by $L(\brho) \triangleq\frac
{1}{1+\veps(\brho)}$,
satisfies $L(\brho) \approx1-\veps(\brho)$ when $\veps= \veps
(\brho)$ is small,
that is, when the system approaches criticality.
Then, an immediate consequence of the bound (\ref{eq:exbdd}) is that
\begin{eqnarray*}
\limsup_{\gamma\rightarrow\infty} \frac{1}{\gamma} \log\PP _{\bpi}\bigl(\|\bN
\|_{\infty} \geq\gamma\bigr) & \lesssim& \frac{1}{2\xi}\log \biggl(
\frac{\xi}{\xi+\veps
K} \biggr)
\\
& \approx& -\frac{K\veps}{2\xi^2} \approx-\frac{K}{2\xi
^2} \bigl(1-L(\brho)
\bigr).
\end{eqnarray*}
Note that $\frac{K}{2\xi^2}$ is a load-dependent constant.
Thus Theorem \ref{thm:sw2} shows that the large-deviations exponent of
the steady-state number of flows is upper bounded by $-(1 - L(\brho))$,
up to a multiplicative constant.

\subsubsection*{Interchange of limits ($\alpha=1$)}
As discussed in the \hyperref[sec:model]{Introduction}, when $\alpha=
1$, Theorem \ref {thm:sw2} leads to the tightness (Lemma
\ref{lem:tight}) of the steady-state distributions of the model under
diffusion scaling. This in turn leads to Theorem \ref{thm:inlimit} and
Corollary~\ref{cor:inlimit}, on the validity of the diffusion
approximation in steady state. As the statements of these results
require a significant amount of preliminary notation and background
(which is {introduced} in
Section \ref{sec:inlimit}), we give {here} an informal statement.

\begin{ilh*}
Consider a sequence of flow-level networks
operating under the proportionally fair policy.
Let $\bN^r(\cdot)$ be the flow-vector Markov process associated with
the $r$th network,
let $\veps^r$ be the corresponding gap, and
let $\hat{\bpi}^r$ be the stationary distribution of $\veps^r\bN
^r(\cdot)$.
As $\veps^r \rightarrow0$, and under certain technical conditions,
$\hat{\bpi}^r$ converges weakly to the stationary distribution of an
associated limiting process.
\end{ilh*}



\section{\texorpdfstring{Transient analysis ($\alpha\geq1$)}
{Transient analysis (alpha >= 1)}}\label{sec4}
In this section, we present a transient analysis of the $\alpha$-fair
policies with $\alpha\geq1$.
First we present a general maximal lemma, which we then specialize to
our model.
In particular, we prove a refined drift inequality for the Lyapunov
function given by
%
\begin{equation}
F_{\alpha}(\mathbf{n}) = \frac{1}{\alpha+1}\sum
_{i\in I}\nu _i\kappa_i
\mu_i^{\alpha-1} \biggl(\frac{n_i}{\nu_i} \biggr)^{\alpha+1}.
\label{flpnv}
\end{equation}
This Lyapunov function and associated drift inequalities have played an
important role in establishing positive recurrence
(cf. \cite{BM01,DeV01,KW04}) and multiplicative
state space collapse (cf.  \cite{KKLW09}) for $\alpha$-fair policies.
We combine our drift inequality with the maximal lemma to obtain a
maximal inequality for bandwidth-sharing networks.
We then apply the maximal inequality to prove full state space collapse
when \mbox{$\alpha\geq1$}.

\subsection{The key lemma}
Our analysis relies on the following lemma.
%
\begin{lemma}\label{maximal}
Let $(\mathscr{F}_n)_{n\in\Zp}$ be a filtration on a probability
space. Let $(X_n)_{n\in\Zp}$ be a nonnegative $\mathscr
{F}_n$-adapted stochastic process that satisfies
%
\begin{equation}
\mathbb{E}[ X_{n+1}\mid\mathscr{F}_n ]\leq
X_n+B_n, \label{eq:maxdrift}
\end{equation}
where the $B_n$ are nonnegative random variables (not necessarily
$\mathscr{F}_n$-adapted) with finite means.
Let $X_n^{*} = \max\{X_0,\ldots,X_n\}$ and suppose that $X_0 = 0$.
Then, for any $a>0$ and any $T\in\Zp$,
\[
\mathbb{P} \bigl(X_T^{*} \geq a \bigr)\leq
\frac{\sum_{n=0}^{T-1}\mathbb{E}[B_n]}{a}.
\]
\end{lemma}
This lemma is a simple consequence of the following standard maximal
inequality for nonnegative supermartingales; see, for example, Exercise
4, Section 12.4, of \cite{grimmett}.
%
\begin{theorem}\label{thm:smart}
Let $(\mathscr{F}_n)_{n\in\Zp}$ be a filtration on a probability
space. Let $(Y_n)_{n\in\Zp}$
be a nonnegative $\mathscr{F}_n$-adapted supermartingale; that is, for
all $n$,
\[
\E[Y_{n+1}\mid\mathscr{F}_n ]\leq Y_n.
\]
Let $Y_T^* = \max\{Y_0,\ldots,Y_T\}$. Then,
\[
\mathbb{P} \bigl(Y_T^*\geq a \bigr)\leq\frac{\E[Y_0]}{a}.
\]
\end{theorem}
\begin{pf*}{Proof of Lemma \ref{maximal}} First note that if we take
the conditional expectation
of both sides of \eqref{eq:maxdrift}, given $\mathscr{F}_n$, we have
\[
\E[X_{n+1}\mid\mathscr{F}_n]\leq\E[X_n\mid
\mathscr{F}_n]+\E [B_n\mid\mathscr{F}_n] =
X_n + \E[B_n\mid\mathscr{F}_n].
\]
Fix $T\in\Zp$. For any $n\leq T$, define
\[
Y_n = X_n + \E \Biggl[\sum
_{k=n}^{T-1}B_k \Big| \mathscr{F}_n
\Biggr].
\]
Then
\begin{eqnarray*}
\E[Y_{n+1}\mid\mathscr{F}_n] & = & \E[X_{n+1}\mid
\mathscr{F}_n] + \E \Biggl[ \E \Biggl[ \sum
_{k=n+1}^{T-1} B_k \Big| \mathscr
{F}_{n+1} \Biggr] \Big| \mathscr{F}_n \Biggr]
\\
&\leq& X_n + \E[B_n\mid\mathscr{F}_n]+ \E
\Biggl[ \sum_{k=n+1}^{T-1} B_k \Big|
\mathscr{F}_n \Biggr] = Y_n.
\end{eqnarray*}
Thus, $Y_n$ is an $\mathscr{F}_n$-adapted supermartingale;
furthermore, by definition, $Y_n$ is nonnegative for all $n$.
Therefore, by Theorem \ref{thm:smart},
\[
\mathbb{P} \bigl(Y_T^*\geq a \bigr)\leq\frac{\E[Y_0]}{a} =
\frac{\E
[\sum_{k=0}^{T-1} B_k ]}{a}.
\]
But $Y_n\geq X_n$ for all $n$, since the $B_k$ are nonnegative. Thus,
\[
\mathbb{P} \bigl(X_T^*\geq a \bigr)\leq\mathbb{P}
\bigl(Y_T^* \geq a \bigr)\leq\frac{\E
[\sum_{k=0}^{T-1} B_k ]}{a}.
\]
\upqed\end{pf*}
Since we are dealing with continuous-time Markov processes, the
following corollary of Lemma \ref{maximal}
will be useful for our analysis.
%
\begin{corollary}\label{cmax}
Let $(\mathscr{F}_t)_{t\geq0}$ be a filtration on a probability space.
Let $Z_t$ be a nonnegative, right-continuous $\mathscr{F}_t$-adapted
stochastic process {that} satisfies
\[
\mathbb{E}[Z_{s+t}|\mathscr{F}_s]\leq Z_s+Bt
\]
for all $s,t\geq0$, where $B$ is a nonnegative constant.
Assume that $Z_0\equiv0$. Denote $Z_T^{*} \triangleq\sup_{0\leq t
\leq T} Z_t$ (which
can possibly be infinite). Then, for any $a>0$, and for any $T\geq0$,
\[
\mathbb{P} \bigl(Z_T^{*} \geq a \bigr)\leq
\frac{BT}{a}.
\]
\end{corollary}
\begin{pf}
The proof is fairly standard. We fix $T\geq0$ and $a>0$.
Since $Z_t$ is right-continuous,
$Z_{T}^{*} = \sup_{t\in[0,T]} Z_t = \sup_{t\in([0,T]\cap\mathbb
{Q})\cup\{T\}} Z_t$.
Consider an increasing sequence of finite sets $I_n$ so that
$\bigcup_{n=1}^{\infty} I_n = ([0,T]\cap\mathbb{Q})\cup\{T\}$,
and $0,T\in I_n$ for all $n$. Define $Z_{T}^{(n)}=\sup_{t\in I_n}Z_t$.
Then $ (Z_{T}^{(n)} )_{n=1}^{\infty}$ is a nondecreasing sequence,
and $Z_T^{(n)}\rightarrow Z_{T}^{*}$ as $n\rightarrow\infty$, almost surely.
For each $Z_{T}^{(n)}$, we can apply Lemma~\ref{maximal},
and it is immediate that for any $b>0$,
%
\begin{equation}
\label{eq:cmax} \mathbb{P} \bigl(Z_{T}^{(n)}> b \bigr)\leq
\frac{BT}{b},
\end{equation}
since each $I_n$ includes both $0$ and $T$.
Since $Z_{T}^{(n)}$ increases monotonically to $Z_{T}^{*}$, almost surely,
we have that $\PP(Z_{T}^{(n)}> b)\leq\PP(Z_{T}^{(n+1)}> b)$ for all $n$,
and $\PP(Z_{T}^{(n)}> b) \rightarrow\PP(Z_{T}^{*}> b)$ as
$n\rightarrow\infty$.
The right-hand side of (\ref{eq:cmax}) is fixed, so
\[
\mathbb{P} \bigl(Z_{T}^{*}> b \bigr)\leq
\frac{BT}{b}.
\]
We now take an increasing sequence $b_n$ with $\lim_{n\rightarrow
\infty} b_n = a$, and obtain
\[
\mathbb{P} \bigl(Z_{T}^{*}\geq a \bigr)\leq
\frac{BT}{a}.
\]
\upqed\end{pf}

\subsection{A maximal inequality for bandwidth-sharing networks}
We employ the Lyapunov function (\ref{flpnv}) to study {$\alpha
$-fair} policies.\vadjust{\goodbreak}
This is the Lyapunov function that was used in \cite{BM01,DeV01} and \cite{KW04}
to establish positive recurrence of the process $\bN(\cdot)$ under
{an} $\alpha$-fair policy.
Below we fine-tune the proof in \cite{DeV01} to obtain a more precise
bound on the Lyapunov drift.
{We note that a ``fluid-model'' version of the following lemma appeared
in the proof of Theorem 1 in
\cite{BM01}.}
For notational convenience, we drop the subscript $\alpha$ from
$F_{\alpha}$
and write $F$ instead.
%
\begin{lemma}\label{lm:max}
Consider a bandwidth-sharing network with $\bA\brho<\bC$
operating under an $\alpha$-fair policy with $\alpha> 0$. Let $\veps
$ be the gap.
Then{,} for any nonzero flow vector $\mathbf{n}$,
\[
\bigl\langle{\nabla F(\mathbf{n}),\bnu-\bmu\bLambda(\mathbf {n})} \bigr\rangle
\leq-\veps \bigl\langle{\nabla F(\mathbf{n}), \bnu} \bigr\rangle,
\]
where $ \langle{\cdot,\cdot} \rangle$ denotes the
standard inner product,
$\nabla F(\mathbf{n})$ denotes the gradient of $F$,
and $\bmu\bLambda(\mathbf{n})$ is the vector $ (\mu_i\Lambda
_i(\mathbf{n}) )_{i\in\cI}$.
\end{lemma}
\begin{pf} We have
\begin{eqnarray*}
\bigl\langle{\nabla F(\mathbf{n}),\bnu-\bmu\bLambda(\mathbf {n})} \bigr
\rangle&=& \sum_{i\in\cI}\frac{1}{\mu_i}\kappa
_i \biggl(\frac{n_i}{\rho_i} \biggr)^{\alpha} \bigl(
\nu_i-\mu_i\Lambda _i(\mathbf{n}) \bigr)
\\
&=& \sum_{i\in I}\kappa_i \biggl(
\frac{n_i}{\rho_i} \biggr)^{\alpha
} \bigl(\rho_i-
\Lambda_i(\mathbf{n}) \bigr)
\\
&=& \bigl\langle{\nabla G_{\mathbf{n}}(\brho_{+}),\brho
_{+}-\bLambda_{+}(\mathbf{n})} \bigr\rangle,
\end{eqnarray*}
where $\brho_{+}=(\rho_i)_{i\in\cI_{+}(\mathbf{n})}$.
Similarly we can get $ \langle{\nabla F(\mathbf{n}),\bnu}
\rangle
=  \langle{\nabla G_{\mathbf{n}}(\brho_{+}),\brho_{+}}
\rangle$.

Now consider the function $g\dvtx [0,1] \rightarrow\mathbb{R}$ defined by
\[
g(\theta)=G_{\mathbf{n}} \bigl(\theta(1+\veps)\brho_{+}+(1-\theta )
\bLambda_{+}(\mathbf{n}) \bigr).
\]
Since $(1+\veps)\brho_{+}$ satisfies the constraints in (\ref
{constraint1}), and
$\bLambda_{+}(\mathbf{n})$ maximizes the strictly \textit{concave}
function $G_{\mathbf{n}}$
subject to the constraints in (\ref{constraint1}), 
we have
\[
G_{\mathbf{n}} \bigl((1+\veps)\brho_{+} \bigr) \leq
G_{\mathbf{n}} \bigl(\bLambda _+(\mathbf{n}) \bigr)\qquad \mbox{that is } g(1)\leq
g(0).
\]
Furthermore, since $G_{\mathbf{n}}$ is a concave function, $g$ is also
concave in $\theta$.
Thus,
\[
g(0) \leq g(1) + (0-1)g'(1) \leq g(0) + (0-1)g'(1).
\]
Hence, $g'(1)\leq0$, that is,
%
\begin{equation}
\label{eq:derv} \frac{dg}{d\theta} \bigg|_{\theta=1} = \bigl\langle{\nabla
G_{\mathbf{n}} \bigl((1+\veps)\brho_{+} \bigr),(1+\veps)
\brho_{+}-\bLambda _{+}(\mathbf{n})} \bigr\rangle\leq0.
\end{equation}
But it is easy to check that $\nabla G_{\mathbf{n}}((1+\veps)\brho
_{+}) = (1+\veps)^{-\alpha}\nabla G_{\mathbf{n}}(\brho_{+})$,
so dividing~(\ref{eq:derv}) by $(1+\veps)^{-\alpha}$, we have
\[
\bigl\langle{\nabla G_{\mathbf{n}}(\brho_{+}),\brho_{+}-
\bLambda _{+}(\mathbf{n})} \bigr\rangle\leq-\veps \bigl\langle{\nabla
G_{\mathbf{n}}(\brho_{+}),\brho_{+}} \bigr\rangle.
\]
This is the same as
\[
\hspace*{82pt}\bigl\langle{\nabla F(\mathbf{n}),\bnu-\bmu\bLambda(\mathbf {n})} \bigr\rangle
\leq-\veps \bigl\langle{\nabla F(\mathbf {n}),\bnu} \bigr\rangle.\hspace*{82pt}\qed\vadjust{\goodbreak}
\]
\noqed\end{pf}
Our next lemma provides a uniform upper bound on the expected
change of $F(\tilde{\bN}(\cdot))$ in one time step,
where $\tilde{\bN}(\cdot)$ is the uniformized chain associated with the
Markov process $\bN(\cdot)$; cf. Definition \ref{df:umc}.
%
\begin{lemma}\label{lm:max2}
Let $\alpha\geq1$. As above, consider a bandwidth-sharing network
with $\bA\brho< \bC$
operating under an $\alpha$-fair policy.
Let $\veps$ be the gap.
Let $ (\tilde{\bN}(\tau) )_{\tau\in\Zp}$ be the
uniformized chain
associated with the Markov process $\bN(\cdot)$.
Then there exists a positive load-dependent constant $\bar{K}$, 
such that for all $\tau\in\Zp$,
\[
\E \bigl[F \bigl(\tilde{\bN}(\tau+1) \bigr)-F(\mathbf{n}) | \tilde{\bN}(\tau ) =
\mathbf{n} \bigr] \leq\bar{K}\veps^{1-\alpha}.
\]
\end{lemma}
\begin{pf}
By the mean value theorem (cf. Proposition \ref{prp:mvt}), for
$\mathbf{n},\bbm\in\mathbb{Z}_{+}^{|\cI|}$, we have
%
\begin{equation}
\label{eq:mvtlypv} F(\mathbf{n}+\bbm)-F(\mathbf{n})= \bigl\langle{\nabla F(
\mathbf {n}),\bbm} \bigr\rangle+\tfrac{1}{2}\bbm^{T}
\nabla^2F(\mathbf {n}+\theta\bbm)\bbm
\end{equation}
for some $\theta\in[0,1]$. We note that, for $\bbm=\pm\be_i$, we have
%
\begin{eqnarray}
\label{eq:2nd_order}
\nonumber
\frac{1}{2}\bbm^{T}
\nabla^2F(\mathbf{n}+\theta\bbm)\bbm 
& \leq&
\frac{\kappa_i\alpha}{2\mu_i\rho_i^{\alpha}}(n_i\pm \theta)^{\alpha-1}
\nonumber
\\[-8pt]
\\[-8pt]
\nonumber
&\leq& \frac{\kappa_i\alpha}{2\mu_i\rho_i^{\alpha
}}(n_i+1)^{\alpha-1},
\end{eqnarray}
since $\alpha\geq1$, and $\theta\in[0,1]$.

As in \cite{DeV01}, we define
\[
\bQ F(\mathbf{n})\triangleq\sum_{\bbm}q(\mathbf{n},
\mathbf {n}+\bbm) \bigl[F(\mathbf{n}+\bbm)-F(\mathbf{n}) \bigr],
\]
so that $\bQ$ is the generator of the Markov process $\bN(\cdot)$.
We now proceed to derive an upper bound for $\bQ F(\mathbf{n})$.
Using equation (\ref{eq:mvtlypv}), we can rewrite $\bQ F(\mathbf{n})$ as
\begin{eqnarray*}
\bQ F(\mathbf{n}) &=& \sum_{\bbm}q(\mathbf{n},
\mathbf{n}+\bbm ) \biggl[ \bigl\langle{\nabla F(\mathbf{n}),\bbm} \bigr\rangle +
\frac{1}{2}\bbm^{T}\nabla^2F(\mathbf{n}+
\theta_{\bbm} \bbm)\bbm \biggr]
\\
&=& \sum_{\bbm}q(\mathbf{n},\mathbf{n}+\bbm) \bigl
\langle{\nabla F(\mathbf{n}),\bbm} \bigr\rangle
\\
& &{} +\frac{1}{2}\sum_{\bbm}q(\mathbf{n},
\mathbf{n}+\bbm )\bbm^{T}\nabla^2F(\mathbf{n}+
\theta_{\bbm} \bbm)\bbm
\end{eqnarray*}
for some scalars $\theta_{\bbm} \in[0,1]$, one such scalar for each
$\bbm$.
From the definition of $\mathbf{q}$, we have
\begin{eqnarray*}
\sum_{\bbm}q(\mathbf{n},\mathbf{n}+\bbm) \bigl
\langle{\nabla F(\mathbf{n}),\bbm} \bigr\rangle & = & \biggl\langle{\nabla F(
\mathbf{n}),\sum_{\bbm}q(\mathbf {n},\mathbf{n}+\bbm)
\bbm} \biggr\rangle
\\
& = & \bigl\langle{\nabla F(\mathbf{n}),\bnu-\bmu\bLambda(\mathbf {n})} \bigr
\rangle.
\end{eqnarray*}
From (\ref{eq:2nd_order}), for $\bbm= \pm\be_i$, we also have
\[
\tfrac{1}{2}\bbm^{T}\nabla^2F(\mathbf{n}+
\theta_{\bbm} \bbm)\bbm \leq\kappa_i\alpha(n_i+1)^{\alpha-1}/2
\mu_i\rho_i^{\alpha}.
\]
Thus
\begin{eqnarray*}
\bQ F(\mathbf{n})&\leq& \bigl\langle{\nabla F(\mathbf{n}),\bnu -\bmu\bLambda(
\mathbf{n})} \bigr\rangle+\sum_{i\in\cI}
\frac
{\kappa_i\alpha}{2\mu_i\rho_i^{\alpha}}(n_i+1)^{\alpha-1} \bigl(\nu _i+
\mu_i\Lambda_i(\mathbf{n}) \bigr)
\\
&\leq& -\veps\sum_{i\in\cI}\kappa_i \biggl(
\frac{n_i}{\rho
_i} \biggr)^{\alpha}\rho_i+\sum
_{i\in\cI}\frac{\kappa_i\alpha
}{2\rho_i^{\alpha}}(n_i+1)^{\alpha-1}
\bigl(\rho_i+\Lambda_i(\mathbf {n}) \bigr)
\\
&\leq& -m\veps\sum_{i\in I}n_i^{\alpha}+M
\sum_{i\in
I}(n_i+1)^{\alpha-1},
\end{eqnarray*}
where the second inequality follows from Lemma \ref{lm:max},
and the third by defining
\[
m\triangleq\min_{i\in\cI}\kappa_i
\rho_i^{1-\alpha}, \qquad M \triangleq\max_{i\in\cI}
\frac{\kappa_i\alpha}{2\rho
_i^{\alpha}} \Bigl(\rho_i+\max_{j\in\cJ}C_j
\Bigr),
\]
{and noting the fact} that since $\Lambda_i(\mathbf{n})\leq\max_{j\in\cJ}C_j$ for all $i$, {we have}
$M \geq\max_{i\in\cI}\frac{\kappa_i\alpha}{2\rho_i^{\alpha
}} (\rho_i+\Lambda_i(\mathbf{n}) )$.
It is then a simple calculation to see that for every $\mathbf{n}\geq
\mathbf{0}$, we have
\[
\bQ F(\mathbf{n}) \leq-m\veps\sum_{i\in I}n_i^{\alpha}+M
\sum_{i\in I}(n_i+1)^{\alpha-1} \leq
\tilde{K}\veps^{1-\alpha}
\]
for some positive load-dependent constant $\tilde{K}$.
Now given $\tilde{\bN}(\tau) = \mathbf{n}$,
\[
\E \bigl[F \bigl(\tilde{\bN}(\tau+1) \bigr)-F(\mathbf{n}) | \tilde{\bN}(\tau ) =
\mathbf{n} \bigr] = \frac{\bQ F(\mathbf{n})}{\Xi} \leq\frac{\tilde{K}\veps^{1-\alpha}}{\Xi}.
\]
%
By setting $\bar{K} = \tilde{K}/\Xi$, we have proved the lemma.
\end{pf}

\begin{corollary}\label{cor:drift}
Let $\alpha\geq1$. As before, suppose that $\bA\brho< \bC$, and
let $\veps$ be the associated gap.
Then, under the $\alpha$-fair policy, the process $\bN(\cdot)$ satisfies
\[
\mathbb{E} \bigl[F \bigl(\bN(s+t) \bigr)-F \bigl(\bN(s) \bigr) | \bN(s) \bigr]
\leq \tilde {K}t\veps^{1-\alpha}\qquad \mbox{for all } t\geq0
\]
for some positive load-dependent constant $\tilde{K}$. 
\end{corollary}
\begin{pf}
The idea of the proof is to show that the expected number of state
transitions of $\bN(\cdot)$
in the time interval $[s, s+t]$ is of order $O(t)$.

Consider the uniformized Markov chain $\tilde{\bN}(\cdot)$
associated with the process $\bN(\cdot)$.
Denote the number of state transitions in the uniformized version of
the process $\bN(\cdot)$
in the time interval $[s, s+t]$ by $\tau$. By the Markov property,
time-homogeneity and the definition of $\tilde{\bN}(\cdot)$, we have
\begin{eqnarray*}
& &\E \bigl[F \bigl(\bN(s+t) \bigr) - F \bigl(\bN(s) \bigr) | \bN(s) =
\mathbf{n} \bigr]
\\
&&\qquad =  \E \bigl[F \bigl(\tilde{\bN}(\tau) \bigr) - F \bigl(\tilde{\bN}(0) \bigr)
| \tilde{\bN}(0) = \mathbf{n} \bigr].
\end{eqnarray*}
Now, by the definition of the uniformized chain, $\tau$ and $\tilde
{\bN}(\cdot)$ are independent. Thus
\begin{eqnarray*}
& &\E \bigl[F \bigl(\tilde{\bN}(\tau) \bigr) - F \bigl(\tilde{\bN}(0) \bigr) |
\tilde{\bN}(0) \bigr]
\\
&&\qquad =  \E \Biggl[\sum_{k=0}^{\tau-1} \bigl(F
\bigl(\tilde{\bN}(k+1) \bigr) - F \bigl(\tilde{\bN}(k) \bigr) \bigr) \Big| \tilde{
\bN}(0) \Biggr]
\\
&&\qquad =  \E \Biggl[\E \Biggl[\sum_{k=0}^{\tau-1}
\bigl(F \bigl(\tilde{\bN }(k+1) \bigr) - F \bigl(\tilde{\bN}(k) \bigr) \bigr) \Big|
\tilde{\bN}(0), \tau \Biggr] \Big| \tilde{\bN}(0) \Biggr]
\\
&&\qquad =  \E \Biggl[\sum_{k=0}^{\tau-1} \E \bigl[F
\bigl(\tilde{\bN}(k+1) \bigr) - F \bigl(\tilde{\bN}(k) \bigr) | \tilde{\bN}(0),
\tau \bigr]\Big | \tilde{\bN}(0) \Biggr]
\\
&&\qquad =  \E \Biggl[\sum_{k=0}^{\tau-1} \E \bigl[F
\bigl(\tilde{\bN}(k+1) \bigr) - F \bigl(\tilde{\bN}(k) \bigr) | \tilde{\bN}(0)
\bigr]\Big | \tilde {\bN}(0) \Biggr]
\\
&&\qquad \leq \E \Biggl[\sum_{k=0}^{\tau-1} \bar{K}
\veps^{1-\alpha
} \Biggr] = \bar{K}\veps^{1-\alpha}\E[\tau]
\end{eqnarray*}
for some load-dependent constant $\bar{K}$.
The fourth equality follows from the independence of $\tau$ and
$\tilde{\bN}(\cdot)$,
and the inequality follows from Lemma \ref{lm:max2}.
Since the counting process of the number of state transitions in the
uniformized version of the process $\bN(\cdot)$
is a time-homogeneous Poisson process of rate $\Xi$, we have $\E[\tau
] = \Xi t$. This shows that
\[
\mathbb{E} \bigl[F \bigl(\bN(s+t) \bigr)-F \bigl(\bN(s) \bigr) | \bN(s) \bigr]
\leq \bar {K}\Xi t\veps^{1-\alpha}.
\]
The proof is completed by setting $\tilde{K} = \bar{K}\Xi$.
\end{pf}

\begin{pf*}{Proof of Theorem \protect\ref{thm:max}}
Let $b>0$. Then
\begin{eqnarray*}
\mathbb{P} \bigl(N^*(T)\geq b \bigr) & = & \mathbb{P} \biggl(\frac
{1}{\alpha+1}
\bigl(N^*(T) \bigr)^{\alpha+1}\geq\frac{1}{\alpha
+1}b^{\alpha+1} \biggr)
\\
&\leq& \mathbb{P} \biggl(\sup_{t\in[0,T]}F \bigl(\bN(t) \bigr)\geq
\biggl(\min_{i\in\cI}\frac{1}{\alpha+1}\kappa_i
\mu_i^{\alpha-1}\nu _i^{-\alpha}
\biggr)b^{\alpha+1} \biggr)
\\
&\leq& \frac{(\alpha+1)K'T}{ (\min_{i\in\cI}\kappa_i \mu
_i^{\alpha-1}\nu_i^{-\alpha} )\veps^{\alpha-1}b^{\alpha+1}} = \frac{KT}{\veps^{\alpha-1}b^{\alpha+1}},
\end{eqnarray*}
where the second inequality follows from Corollaries \ref{cmax} and
\ref{cor:drift},
$K'$ is as in Corollary \ref{cor:drift},
and $K = \frac{(\alpha+1)K'}{\min_{i\in\cI}\kappa_i \mu
_i^{\alpha-1}\nu_i^{-\alpha}}$.
\end{pf*}

\subsection{\texorpdfstring{Full state space collapse for $\alpha\geq1$}
{Full state space collapse for alpha >= 1}}\label{ssec:ssc}
Throughout this section, we assume that we have fixed $\alpha\geq1$,\vadjust{\goodbreak}
and correspondingly, the Lyapunov function (\ref{flpnv}).
To state the full state space collapse result for $\alpha\geq1$,
we need some preliminary definitions and
the statement of the multiplicative state space collapse result.

Consider a sequence of bandwidth-sharing networks indexed by $r$,
where $r$ is to be thought of as increasing to infinity along a sequence.
Suppose that the incidence matrix $\bA$, the capacity vector $\bC$
and the weights
$\{\kappa_i\dvtx i\in\cI\}$ do not vary with~$r$. Write $\bN^{r}(t)$ for
the flow-vector
Markov process associated with the $r$th network. Similarly, we write
$\bnu^r$, $\bmu^r$, $\brho^r$, etc.
We assume the following \textit{heavy-traffic} condition (cf. \cite{KKLW09}):
%
\begin{asmption}\label{asmp:fssc}
We assume that $\bA\brho^r < \bC$ for all $r$. We also assume that
there exist $\bnu,\bmu\in\mathbb{R}_+^{|\cI|}$ and $\btheta>
\mathbf{0}$,
such that $\nu_i>0$ and $\mu_i>0$ for all $i\in\cI$, $\bnu
^r\rightarrow\bnu$ and
$\bmu^r\rightarrow\bmu$ as $r\rightarrow\infty$, and $r(\bC-\bA
\brho^r)\rightarrow\btheta$ as $r\rightarrow\infty$.
\end{asmption}
Note that our assumption differs from that in \cite{KKLW09}, which
allows convergence to the critical load
from both overload and underload, whereas here we only allow
convergence to the critical load from underload.

To state the multiplicative state space collapse result,
we also need to define a \textit{workload process} $\bW(t)$ and a \emph
{lifting map} $\Delta$.
%
\begin{definition}\label{df:lift}
We first define the \textit{workload} $\mathbf{w}\dvtx \Rpp^{|\cI|}
\rightarrow\Rpp^{|\cJ|}$
associated with a flow-vector $\mathbf{n}$ by
$\mathbf{w}= \mathbf{w}(\mathbf{n}) = \bA\bM^{-1}\mathbf{n}$,
where $\bM= \diag(\bmu)$ is the $|\cI|\times|\cI|$ diagonal
matrix with $\bmu$ on its diagonal.
The \textit{workload process} $\bW(t)$ is defined to be $\bW(t)
\triangleq\bA\bM^{-1}\bN(t)$, for all $t\geq0$.
We also define the \textit{lifting map} $\Delta$. For each $\mathbf
{w}\in\mathbb{R}_+^{|\cJ|}$, define
$\Delta(\mathbf{w})$ to be the unique value of $\mathbf{n}\in
\mathbb{R}_+^{|\cI|}$ that solves the following optimization problem:
\begin{eqnarray*}
&&\mathrm{minimize}\  F(\mathbf{n}) \\
&&\mathrm{subject\ to} \qquad \sum_{i\in\cI}A_{ji}\frac{n_i}{\mu_i}\geq
w_j,  j\in\cJ,  n_i\geq0,  i\in\cI.
\end{eqnarray*}
\end{definition}
For simplicity, suppose that all networks start with zero flows.
We consider the following diffusion scaling:
%
\begin{equation}
\label{scaling} \hat{\bN}^r(t)=\frac{\bN^r(r^2t)}{r} \quad\mathrm{and}\quad
\hat{\bW }^r(t)=\frac{\bW^r(r^2t)}{r},
\end{equation}
where $\bW^r(t) = \bA(\bM^r)^{-1}{\bN^r(t)}$, and $\bM^r = \diag
(\bmu^r)$.

The following multiplicative state space collapse result is known to hold.
%
\begin{theorem}[(Multiplicative state space collapse \cite{KKLW09}, Theorem
5.1)]\label{thm:mssc}
Fix $T>0$ and assume that $\alpha\geq1$.
Write $\|\mathbf{x}(\cdot)\|=\sup_{t\in[0,T],i \in\cI}|x_i(t)|$.
Then, under Assumption~\ref{asmp:fssc},
and for any $\delta>0$,
\[
\lim_{r\rightarrow\infty} \mathbb{P} \biggl( \frac{\|\hat{\bN
}^r(\cdot)-\Delta(\hat{\bW}^r(\cdot))\|}{\|\hat{\bN}^r(\cdot)\|
}>\delta
\biggr)=0.\vadjust{\goodbreak}
\]
\end{theorem}
We can now state and prove a full state space collapse result:
%
\begin{theorem}[(Full state space collapse)]\label{thm:fssc}
Under the same assumptions as in Theorem \ref{thm:mssc}, and for any
$\delta>0$,
\[
\lim_{r\rightarrow\infty} \mathbb{P} \bigl( \bigl\|\hat{\bN}^r(\cdot
)-\Delta \bigl(\hat{\bW}^r(\cdot) \bigr)\bigr\|>\delta \bigr)= 0.
\]
\end{theorem}
\begin{pf}
Let $\veps_r = \veps(\brho^r)$ be the gap in the $r$th system. Then,
under Assumption~\ref{asmp:fssc},
$\veps_r\geq D/r$ for some network-dependent constant $D>0$, and for
$r$ sufficiently large.
By Theorem \ref{thm:max}, for any $b>0$, and for sufficiently large $r$,
\begin{eqnarray*}
\mathbb{P} \bigl( N^{r,*} \bigl(r^2T \bigr)\geq b \bigr) &
\leq& \frac{K_r
r^2T}{\veps_r^{\alpha-1}b^{\alpha+1}}
\\
&\leq& \frac{K_r r^{1+\alpha} T}{D^{\alpha-1} b^{\alpha+1}}.
\end{eqnarray*}
Here, $K_r$ is a load-dependent constant associated with the $r$th system,
as specified in the proof of Theorem \ref{thm:max}.
From the proof of Theorem \ref{thm:max}, note also that
$K_r = f(\bmu^r, \bnu^r)$, for a function $f$
that is continuous on the open positive orthant $\R_p^{|\cI|}\times
\R_p^{|\cI|}$.
Since $\bmu^r \rightarrow\bmu> \mathbf{0}$, and $\bnu^r
\rightarrow\bnu> \mathbf{0}$,
$K_r \rightarrow K \triangleq f(\bmu, \bnu) \in\R$.
In particular, the $K_r$ are bounded, and for all sufficiently large~$r$,
\[
\mathbb{P} \bigl( N^{r,*} \bigl(r^2T \bigr)\geq b \bigr)
\leq \frac{ (K+1)
r^{1+\alpha} T}{D^{\alpha-1} b^{\alpha+1}}.
\]
Then, with $a=b/r$ and under the scaling in (\ref{scaling}),
%
\begin{equation}
\label{eq:excursion} \mathbb{P} \bigl( \bigl\| \hat{\mathbf{N}}^r(\cdot)
\bigr\| \geq a \bigr)\leq \frac{K+1}{D^{\alpha-1}} \cdot\frac{T}{a^{\alpha+1}}
\end{equation}
for any $a>0$.

For notational convenience, we write
\[
B(r) = \bigl\|\hat{\mathbf{N}}^r(\cdot)-\Delta \bigl(\bW^{r}(
\cdot) \bigr)\bigr\|.
\]
Then, for any $a>1$,
and for sufficiently large $r$,
\begin{eqnarray*}
\mathbb{P} \bigl( B(r) > \delta \bigr)&\leq& \mathbb{P} \biggl(
\frac
{B(r)}{\|\hat{\mathbf{N}}^r(\cdot)\|}> \frac{\delta}{a} \mbox{ or }  \bigl\| \hat{\mathbf{N}}^r(
\cdot)\bigr\| \geq a \biggr)
\\
&\leq& \mathbb{P} \biggl( \frac{B(r)}{\|\hat{\mathbf{N}}^r(\cdot)\|
}> \frac{\delta}{a} \biggr)+
\mathbb{P} \bigl(\bigl\|\hat{\mathbf {N}}^r(\cdot)\bigr\| \geq a \bigr).
\end{eqnarray*}
Note that by Theorem \ref{thm:mssc}, the first term on the right-hand
side goes to $0$ as $r\rightarrow\infty$, for any $a>0$.
The second term on the right-hand side can be made smaller than any,
arbitrarily small, constant (uniformly, for all $r$), by taking $a$
sufficiently large; cf. equation (\ref{eq:excursion}).
Thus, $\mathbb{P}(B(r)\geq\delta)\rightarrow0$ as $r\rightarrow
\infty$. This completes the proof.
\end{pf}

\section{\texorpdfstring{$\alpha$-fair policies: A useful drift inequality}
{alpha-fair policies: A useful drift inequality}}\label{sec:drift}

We now {shift} 
our focus to the steady-state regime.
The key to many of our results 
is a \textit{drift inequality} that holds
for every $\alpha>0$ and every $\brho> \mathbf{0}$ with $\bA\brho<
\bC$.
In this section,
we shall state and prove this inequality. It will be used in Section
\ref{sec:exp} to
prove Theorem \ref{thm:sw2}.

We define the Lyapunov function that we will employ.
For $\alpha\geq1$, it will be
simply the weighted $(\alpha+1)$-norm $L_{\alpha}(\mathbf{n}) =
\sqrt[\alpha+1]{(\alpha+1)F_{\alpha}(\mathbf{n})}$
of a vector $\mathbf{n}$, where $F_{\alpha}$ was defined in (\ref{flpnv}).
However, when $\alpha\in(0,1)$, this function has unbounded second
derivatives as we
approach the boundary of $\mathbb{R}_+^{|\cI|}$. For this reason, our
Lyapunov function will
be a suitably smoothed version of $\sqrt[\alpha+1]{(\alpha
+1)F_{\alpha}(\cdot)}$.

\begin{definition}\label{df:sw-lyap}
Define $h_{\alpha} \dvtx \mathbb{R}_+ \rightarrow\mathbb{R}_+$
to be $h_{\alpha}(r) = r^{\alpha}$, when $\alpha\geq1$, and
\[
h_{\alpha}(r) = \cases{ %
r^{\alpha},
&\quad $\mbox{if }  r\geq1,$
\vspace*{2pt}\cr
(\alpha-1) r^3 +(1-\alpha)r^2 + r, & \quad $\mbox{if }  r <
1,$}
\]
when $\alpha\in(0,1)$. Let $H_{\alpha} \dvtx \mathbb{R}_+ \rightarrow
\mathbb{R}_+$
be the antiderivative of $h_{\alpha}$, so that $H_{\alpha}(r) = \int_{0}^{r} h_{\alpha}(s)  \,ds$.
The Lyapunov function $L_{\alpha} \dvtx \mathbb{R}_+^{|\cI|} \rightarrow
\mathbb{R}_+$ is defined
to be
\[
L_{\alpha} (\mathbf{n}) = \biggl[ (\alpha+1)\sum
_{i\in\cI} \kappa _i \mu_i^{\alpha-1}
\nu_i^{-\alpha}H_{\alpha}(n_i)
\biggr]^{
{1}/{(\alpha+1)}}.
\]
\end{definition}
For notational convenience, define
%
\begin{equation}
\label{defn:weight} w_i = \kappa_i
\mu_i^{\alpha-1} \nu_i^{-\alpha} \qquad\mbox{for
each } i\in\cI,
\end{equation}
so that more compactly, we have
\[
F_{\alpha}(\mathbf{n}) = \frac{1}{\alpha+1}\sum
_{i\in\cI}w_i n_i^{\alpha+1}\quad \mbox{and}\quad L_{\alpha}(\mathbf{n}) = \biggl[(\alpha+1)\sum
_{i\in
\cI}w_i H_{\alpha}(n_i)
\biggr]^{1/(\alpha+1)}.
\]

We will make heavy use of various properties
of the functions $h_{\alpha}$, $H_{\alpha}$ and $L_{\alpha}$, which
we summarize in the following lemma. The proof is elementary and is omitted.

\begin{lemma}\label{lm:prpt}
Let $\alpha\in(0,1)$. The function $h_{\alpha}$ has the following properties:
\begin{longlist}[(iii)]
\item[(i)] it is continuously differentiable with
$h_{\alpha}(0)=0, h_{\alpha}(1)=1, h'_{\alpha}(0)=1$ and
$h'_{\alpha}(1) = \alpha$;
\item[(ii)] it is increasing and, in particular, $h_{\alpha}(r)\geq
0$ for all $r\geq0$;
\item[(iii)] we have $r^{\alpha}-1\leq h_{\alpha}(r)\leq r^{\alpha
}+1$, for all $r\in[0,1]$;
\item[(iv)] $h'_{\alpha}(r)\leq2$, for all $r\geq0$.\vadjust{\goodbreak}
\end{longlist}
Furthermore, from \textup{(iii)}, we also have the following property of
$H_{\alpha}$:
\begin{longlist}[(iii$'$)]
\item[(iii$'$)] $r^{\alpha+1}-2\leq(\alpha+1)H_{\alpha}(r)\leq
r^{\alpha+1}+2$ for all $r\geq0$.\vspace*{-3pt}
\end{longlist}
\end{lemma}

We are now ready to state the drift inequality.
Here we consider the uniformized chain $ (\tilde{\bN}(\tau
) )_{\tau\in\Zp}$
associated with $\bN(\cdot)$, and the corresponding drift.\vspace*{-3pt}
%
\begin{theorem}\label{thm:sw-drift}
Consider a bandwidth-sharing network operating under an $\alpha$-fair
policy with $\alpha> 0$, and assume that
$\bA\brho< \bC$. Let $\veps$ be the gap.
Then, there exists a positive constant $B$ and a positive
load-dependent constant
$K$, such that if $L_{\alpha}(\tilde{\bN}(\tau))>B$, then
%
\begin{equation}
\label{eq:sw-drift} \mathbb{E} \bigl[L_{\alpha} \bigl(\tilde{\bN}(\tau+1)
\bigr) - L_{\alpha} \bigl(\tilde {\bN}(\tau) \bigr) \mid\tilde{\bN}(\tau)
\bigr] \leq-\veps K.
\end{equation}
Furthermore, $B$ takes the form $K'/\veps$ when $\alpha\geq1$, and
$K'/\min\{\veps^{1/\alpha}, \veps\}$
when $\alpha\in(0,1)$, with $K'$ being a positive load-dependent constant.\vspace*{-3pt}
\end{theorem}
As there is a marked difference between the form of $L_{\alpha}$ for
the two cases $\alpha\geq1$
and $\alpha\in(0,1)$, the proof of the drift inequality is split into
two parts.
We first prove the drift inequality when $\alpha\geq1$, in which case
$L_{\alpha}$ takes a nicer form,
and we can apply results on $F_{\alpha}$ from previous sections. The
proof for the case $\alpha\in(0,1)$
is similar but more tedious.
{We note that such a qualitative difference between the two cases,
$\alpha< 1$ and $\alpha\geq1$,
has also been observed in other works, such as, for example, \cite{Srikant04}.}

We wish to draw attention here to the main difference from related
drift inequalities in the literature.
The usual proof of stability involves the Lyapunov function~(\ref{flpnv});
for instance, for the $\alpha$-fair policy with $\alpha=1$ (the
proportionally fair policy),
it involves a weighted quadratic Lyapunov function. In contrast, we use
$L_{\alpha}$,
a weighted norm function (or its smoothed version),
which scales linearly along radial directions. In this sense, our
approach is
similar in spirit to \cite{BGT01}, which employed piecewise linear
Lyapunov functions
to derive drift inequalities and then moment and tail bounds.
{The use of normed Lyapunov functions to establish stability and
performance bounds
has also been considered in other works; see, for example, \cite{VL09}
and \cite{ES11}.}\vspace*{-3pt}

\subsection{\texorpdfstring{Proof of Theorem \protect\ref{thm:sw-drift}: $\alpha\geq1$}
{Proof of Theorem 5.3: alpha >= 1}}

We wish to decompose the drift term in (\ref{eq:sw-drift})
into the sum of a first-order term and a second-order term, and
we accomplish this by using the second-order mean value theorem; cf. Proposition \ref{prp:mvt}.
Throughout this proof, we drop the subscript $\alpha$
from $L_{\alpha}$ and $F_{\alpha}$, and write $L$ and $F$, respectively.

Consider the function
$L(\mathbf{n})= (\sum_{i\in\cI} w_i n_i^{\alpha+1}
)^{{1}/{(\alpha+1)}}
= [(\alpha+1)F(\mathbf{n}) ]^{{1}/{(\alpha+1)}}.$
The first derivative of $L$ with respect to $\mathbf{n}$ is
$\nabla L(\mathbf{n}) = \nabla F(\mathbf{n})/L^{\alpha}(\mathbf{n})$
by the chain rule and the definition of $L$.
The second derivative is
\[
\nabla^2 L(\mathbf{n}) = \frac{\nabla^2 F(\mathbf{n})}{L^{\alpha
}(\mathbf{n})} - \frac{\nabla F(\mathbf{n}) \nabla L^{\alpha
}(\mathbf{n})^T}{L^{2\alpha}(\mathbf{n})}
= \frac{\nabla^2 F(\mathbf{n})}{L^{\alpha}(\mathbf{n})} - \alpha \frac{\nabla F(\mathbf{n})\nabla F(\mathbf{n})^T}{L^{2\alpha
+1}(\mathbf{n})},
\]
by the quotient rule and the chain rule.\vadjust{\goodbreak}

Write $\mathbf{n}$ for $\tilde{\bN}(\tau)$ and $\mathbf{n}+\bbm$
for $\tilde{\bN}(\tau+1)$,
so that $\bbm= \tilde{\bN}(\tau+1) - \tilde{\bN}(\tau)$.
By Proposition \ref{prp:mvt}, for some $\theta\in[0,1]$, we have
%
\begin{eqnarray}
L(\mathbf{n}+\bbm)-L(\mathbf{n}) &=& \bbm^T \nabla L(\mathbf{n}) +
\frac{1}{2}\bbm^T \nabla^2 L(\mathbf{n}+\theta\bbm)
\bbm
\\
&=& \frac{\bbm^T \nabla F(\mathbf{n})}{L^{\alpha}(\mathbf{n})} + \frac{1}{2} \frac{\bbm^T \nabla^2 F(\mathbf{n}+\theta\bbm) \bbm
}{L^{\alpha}(\mathbf{n}+\theta\bbm)}
\\
& & {}- \frac{\alpha}{2} \frac{\bbm^T\nabla F(\mathbf{n}+\theta\bbm
)\nabla F(\mathbf{n}+\theta\bbm)^T\bbm}{L^{2\alpha+1}(\mathbf
{n}+\theta\bbm)}
\\
&\leq& \frac{\bbm^T \nabla F(\mathbf{n})}{L^{\alpha}(\mathbf{n})} + \frac{1}{2} \bbm^T
\frac{\nabla^2 F(\mathbf{n}+\theta\bbm
)}{L^{\alpha}(\mathbf{n}+\theta\bbm)} \bbm, \label{eq:mvt1}
\end{eqnarray}
since the term $\bbm^T\nabla F(\mathbf{n}+\theta\bbm)\nabla
F(\mathbf{n}+\theta\bbm)^T\bbm$ is nonnegative.
We now consider the two terms in (\ref{eq:mvt1}) separately. Recall
from the proof of Lemma \ref{lm:max2} that
\[
\E \bigl[\bbm^T \nabla F(\mathbf{n}) | \mathbf{n} \bigr] =
\frac{ \langle{\nabla F(\mathbf{n}),\bnu-\bmu\bLambda
(\mathbf{n})} \rangle}{\Xi} \leq-\veps\frac{ \langle{\nabla F(\mathbf{n}),\bnu}
\rangle}{\Xi}.
\]
But $ \langle{\nabla F(\mathbf{n}),\bnu} \rangle= \sum_{i\in\cI}w_i\nu_i n_i^{\alpha}$,
so
%
\begin{equation}
\label{eq:1order} \E \bigl[\bbm^T \nabla F(\mathbf{n}) | \mathbf{n}
\bigr] \leq -\veps\frac{\sum_{i\in\cI}w_i\nu_i n_i^{\alpha}}{\Xi},
\end{equation}
and so
%
\begin{eqnarray}\label{eq:drift1_geq1}
\E \biggl[\frac{\bbm^T \nabla F(\mathbf{n})}{L^{\alpha}(\mathbf
{n})} \Big| \mathbf{n} \biggr] &\leq& -\veps
\frac{\sum_{i\in\cI}w_i\nu_i n_i^{\alpha}}{\Xi
(\sum_{i\in\cI}w_i n_i^{\alpha+1} )^{{\alpha
}/{(\alpha+1)}}}
\nonumber
\\
& = & -\veps\frac{\sum_{i\in\cI}w_i\nu_i n_i^{\alpha}}{\Xi
(\sum_{i\in\cI} (w_i^{{1}/{(\alpha+1)}} n_i )^{\alpha
+1} )^{{\alpha}/{(\alpha+1)}}}
\nonumber
\\
&\leq& -\veps\frac{\sum_{i\in\cI}w_i\nu_i n_i^{\alpha}}{\Xi
\cdot\sum_{i\in I}w_i^{{\alpha}/{(\alpha+1)}} n_i^{\alpha}}
\nonumber
\\[-8pt]
\\[-8pt]
\nonumber
&\leq& -\veps\frac{\max_{i\in\cI}w_i^{{1}/{(\alpha+1)}}\nu
_i}{\Xi}
\\
&=& -\veps\frac{\max_{i\in\cI}\kappa^{{1}/{(\alpha+1)}}\mu
_i^{{(\alpha-1)}/{(\alpha+1)}}\nu_i^{{1}/{(\alpha+1)}}}{\Xi}
\nonumber
\\
&=& -\veps K,\nonumber
\end{eqnarray}
where
%
\begin{equation}
\label{eq:Kconst} { K = K(\alpha, \bkappa, \bmu, \bnu) \triangleq
\frac{\max_{i\in
\cI}\kappa^{{1}/{(\alpha+1)}}\mu_i^{{(\alpha-1)}/{(\alpha
+1)}}\nu_i^{{1}/{(\alpha+1)}}}{\Xi} }
\end{equation}
is a positive load-dependent constant.
The second inequality follows from the fact that for any vector
$\mathbf{x}$, and for any $\alpha> 0$,
$\|\mathbf{x}\|_{\alpha+1}\leq\|\mathbf{x}\|_{\alpha}$.
The second to last equality follows from the definition of the $w_i$;
cf. equation (\ref{defn:weight}).

For the second term in (\ref{eq:mvt1}), we wish to show that
if $L(\mathbf{n})$ is sufficiently large, then
\[
\frac{1}{2} \bbm^T \frac{\nabla^2 F(\mathbf{n}+\theta\bbm
)}{L^{\alpha}(\mathbf{n}+\theta\bbm)} \bbm \leq
\frac{\veps}{2} K.
\]
Note that with probability $1$,
either $\bbm= \mathbf{0}$ or $\bbm= \pm\be_i$ for some $i\in
\cI$.
Thus
\begin{eqnarray*}
\frac{1}{2} \bbm^T \frac{\nabla^2 F(\mathbf{n}+\theta\bbm
)}{L^{\alpha}(\mathbf{n}+\theta\bbm)} \bbm &\leq&
\frac{1}{2}\frac{\max_{i\in\cI} [\nabla^2 F(\mathbf
{n}+\theta\bbm) ]_{ii}}{L^{\alpha}(\mathbf{n}+\theta\bbm)}
\\
&=& \frac{\alpha}{2}\frac{\max_{i\in\cI}w_i (n_i+\theta
m_i)^{\alpha-1}}{ [\sum_{i\in\cI}w_i (n_i+\theta m_i)^{\alpha
+1}  ]^{{\alpha}/{(\alpha+1)}}}
\\
&\leq& \frac{\alpha}{2}\frac{\max_{i\in\cI} w_i(n_i+\theta
m_i)^{\alpha-1}}{w_{i_0}^{{\alpha}/{(\alpha+1)}} (n_{i_0}+\theta
m_{i_0})^{\alpha}}
\\
&\leq& \frac{\alpha}{2}w_{i_0}^{{1}/{(\alpha+1)}} (n_{i_0}+
\theta m_{i_0})^{-1}
\\
&\leq& \frac{\alpha}{2}\max_{i\in\cI}w_i^{{1}/{(\alpha+1)}}
(n_{i_0}+\theta m_{i_0})^{-1},
\end{eqnarray*}
where $i_0 \in\cI$ is such that $ w_{i_0}(n_{i_0}+\theta
m_{i_0})^{\alpha-1} = \max_{i\in\cI} w_i(n_i+\theta m_i)^{\alpha-1}$.

{Now note that
\[
\frac{\alpha}{2}\max_{i\in\cI}w_i^{{1}/{(\alpha+1)}}
(n_{i_0}+\theta m_{i_0})^{-1} \leq
\frac{\veps}{2}K
\]
[where $K$ is defined in (\ref{eq:Kconst})] if and only if
\[
n_{i_0}+\theta m_{i_0} \geq\frac{\alpha\max_{i\in\cI}w_i^{
{1}/{(\alpha+1)}}}{K}\cdot
\frac{1}{\veps},
\]
which holds if $L(\mathbf{n}) \geq K'/\veps$ for some appropriately
defined load-dependent constant $K'$.
Thus, if $L(\mathbf{n}) \geq K'/\veps$, then
%
\begin{equation}
\label{eq:drift2_geq1} \frac{1}{2} \bbm^T
\frac{\nabla^2 F(\mathbf{n}+\theta\bbm
)}{L^{\alpha}(\mathbf{n}+\theta\bbm)} \bbm \leq\frac{\veps}{2}K.
\end{equation}
By adding (\ref{eq:drift1_geq1}) and (\ref{eq:drift2_geq1}),
we conclude that
\[
\E \bigl[L(\mathbf{n}+\bbm) - L(\mathbf{n}) | \mathbf{n} \bigr] \leq-
\frac{\veps}{2} K,
\]
when $L(\mathbf{n}) \geq K'/\veps$.} \qed


\subsection{\texorpdfstring{Proof of Theorem \protect\ref{thm:sw-drift}: $\alpha\in(0,1)$}
{Proof of Theorem 5.3: alpha in (0,1)}}
The proof in this section is similar to that for the case $\alpha\geq1$.
We invoke Proposition \ref{prp:mvt} to write the drift term as a sum
of terms, which we bound separately. As in the previous section,
we drop the subscript $\alpha$ from $L_{\alpha}$, $F_{\alpha}$,
$H_{\alpha}$ and $h_{\alpha}$,
and write instead $L$, $F$, $H$, and $h$, respectively.
Note that to use Proposition \ref{prp:mvt},
we need $L$ to be twice continuously differentiable. Indeed, by Lemma
\ref{lm:prpt}(i),
$h$ is continuously differentiable, so its antiderivative $H$ is twice
continuously
differentiable, and so is $L$. Thus, by the second order mean value theorem,
we obtain an equation similar to equation (\ref{eq:mvt1}),
%
\begin{eqnarray}
L(\mathbf{n}+\bbm) - L(\mathbf{n}) &=& \bbm^T \nabla L(\mathbf{n}) +
\frac{1}{2}\bbm^T \nabla^2 L(\mathbf{n}+\theta
\bbm) \bbm
\\
&\leq& \frac{\sum_{i\in I} m_i w_i h(n_i)}{L^{\alpha}(\mathbf{n})} + \frac{1}{2}\frac{\sum_{i\in\cI} m_i^2 w_i h'(n_i+\theta
m_i)}{L^{\alpha}(\mathbf{n}+\theta\bbm)}
\\
&\leq& \frac{\sum_{i\in I} m_i w_i h(n_i)}{L^{\alpha}(\mathbf{n})} + \frac{1}{2}\frac{\max_{i\in\cI} w_i h'(n_i+\theta m_i)}{L^{\alpha
}(\mathbf{n}+\theta\bbm)} \label{eq:mvt2}
\end{eqnarray}
for some constant $\theta\in[0,1]$, and
where, as before, $\tilde{\bN}(\tau) = \mathbf{n}$ and $\tilde{\bN
}(\tau+1) = \mathbf{n}+\bbm$,
and the last inequality follows from the fact that with probability 1,
either $\bbm= \mathbf{0}$, or $\bbm= \pm\be_i$, for some $i\in
\cI$, and that $h'$
is nonnegative.

We now bound the two terms in (\ref{eq:mvt2}) separately.
{Let us first concentrate on the term
\[
\frac{\sum_{i\in I} m_i w_i h(n_i)}{L^{\alpha}(\mathbf{n})}.
\]
}
By Lemma \ref{lm:prpt}(iii),
\[
\sum_{i\in I} m_i w_i
h(n_i) \leq\sum_{i\in I} m_i
w_i \bigl(n_i^{\alpha}+1 \bigr) \leq\sum
_{i\in I} m_i w_i n_i^{\alpha}
+ \sum_{i\in I} m_i w_i,
\]
{so
\[
\frac{\sum_{i\in I} m_i w_i h(n_i)}{L^{\alpha}(\mathbf{n})} \leq\frac{\sum_{i\in I} m_i w_i n_i^{\alpha}}{L^{\alpha}(\mathbf{n})} + \frac{\sum_{i\in I} m_i w_i}{L^{\alpha}(\mathbf{n})}.
\]
First consider the term $\frac{\sum_{i\in I} m_i w_i n_i^{\alpha
}}{L^{\alpha}(\mathbf{n})}$.}
Note that $\sum_{i\in I} m_i w_i n_i^{\alpha} = \bbm^T \nabla
F(\mathbf{n})$.
We also recall from the proof of Lemma \ref{lm:max} that
\[
\E \bigl[\bbm^T \nabla F(\mathbf{n}) | \mathbf{n} \bigr] =
\frac{ \langle{\nabla F(\mathbf{n}),\bnu-\bmu\bLambda
(\mathbf{n})} \rangle}{\Xi} \leq-\veps\frac{ \langle{\nabla F(\mathbf{n}),\bnu}
\rangle}{\Xi}.
\]
We then proceed along the same lines as in the case $\alpha\geq1$,
{and obtain
that if $L(\mathbf{n}) \geq K_2/\veps$ for some positive
load-dependent constant $K_2$, then}
%
\begin{eqnarray}\label{eq:drift1_leq1}
&&\E \biggl[\frac{\sum_{i\in I} m_i w_i n_i^{\alpha}}{L^{\alpha
}(\mathbf{n})} \Big| \mathbf{n} \biggr]\nonumber\\
&&\qquad\leq-\frac{3}{4}
\veps\frac{\max_{i\in\cI}w_i^{
{1}/{(\alpha+1)}}\nu_i}{\Xi}
\\
&&\qquad= -\frac{3}{4}\veps\frac{\max_{i\in\cI}\kappa^{
{1}/{(\alpha+1)}}\mu_i^{{(\alpha-1)}/{(\alpha+1)}}\nu_i^{
{1}/{(\alpha+1)}}}{\Xi}\nonumber\\
&&\qquad= -\frac{3}{4}\veps K.\nonumber
\end{eqnarray}
Here as in the proof for the case $\alpha\geq1$,
\[
K = K(\alpha,
\bkappa, \bmu, \bnu) \triangleq\frac{\max_{i\in\cI}\kappa
^{{1}/{(\alpha+1)}}\mu_i^{{(\alpha-1)}/{(\alpha+1)}}\nu
_i^{{1}/{(\alpha+1)}}}{\sum_{i\in\cI} \nu_i}
\]
is a positive load-dependent constant.

{Now consider the term $\frac{\sum_{i\in I} m_i w_i}{L^{\alpha
}(\mathbf{n})}$.}
With probability $1$, either $\bbm= \mathbf{0}$ or $\bbm= \pm
\be_i$ for some $i\in\cI$,
and therefore $\sum_{i\in\cI}m_i w_i \leq\max_{i\in\cI}w_i$. Thus,
\[
\E \biggl[\frac{\sum_{i\in I} m_i w_i h(n_i)}{L^{\alpha}(\mathbf
{n})} \Big| \mathbf{n} \biggr] \leq-\frac{3}{4}\veps K
+ \frac{\max_{i\in\cI}w_i}{L^{\alpha
}(\mathbf{n})}.
\]

For the second term in (\ref{eq:mvt2}), note that with $\alpha\in(0,1)$,
Lemma \ref{lm:prpt}(iv) implies that
$h' \leq2$, and therefore,
\[
\frac{1}{2}\frac{\max_{i\in\cI} w_i h'(n_i+\theta m_i)}{L^{\alpha
}(\mathbf{n}+\theta\bbm)} \leq\frac{\max_{i\in\cI} w_i}{L^{\alpha}(\mathbf{n}+\theta\bbm)}.
\]
Note that $L^{\alpha}(\mathbf{n}+\theta\bbm)$ and $L^{\alpha
}(\mathbf{n})$ differ only by a load-dependent constant,
since with probability $1$, either $\bbm= \mathbf{0}$ or $\bbm=
\pm\be_i$ for some $i\in\cI$.
Thus, if $L^{\alpha}(\mathbf{n})\geq K_3/\veps$ for some positive
load-dependent constant $K_3$, then
%
\begin{equation}
\label{eq:drift2_leq1} \frac{\max_{i\in\cI}w_i}{L^{\alpha}(\mathbf{n})} + \frac{\max_{i\in\cI} w_i}{L^{\alpha}(\mathbf{n}+\theta\bbm)} \leq
\frac{1}{4}\veps K.
\end{equation}
Putting (\ref{eq:drift1_leq1}) and (\ref{eq:drift2_leq1}) together,
we get that if $L(\mathbf{n}) \geq K'/\min\{\veps^{1/\alpha},\veps
\}$, where $K' = \max\{K_3^{1/\alpha}, K_2\}$, then
\[
\E \bigl[L(\mathbf{n}+\bbm) - L(\mathbf{n}) | \mathbf{n} \bigr] \leq-
\frac{\veps}{2} K. \qed
\]
%

\section{\texorpdfstring{Exponential tail bound under $\alpha$-fair policies}
{Exponential tail bound under alpha-fair policies}}\label{sec:exp}
In this section, we derive an exponential upper bound on the tail probability
of the stationary distribution of the flow sizes, under an $\alpha
$-fair policy with $\alpha> 0$.
We will use the following theorem, a modification of Theorem 1 from
\cite{BGT01}.
%
\begin{theorem}\label{thm:bgt01}
Let $\mathbf{X}(\cdot)$ be an {irreducible and aperiodic}
discrete-time Markov chain with a countable state space $\mathscr{X}$.
Suppose that there exists a Lyapunov function $\Phi\dvtx  \mathscr{X}
\rightarrow\mathbb{R}_+$ with the following properties:
\begin{longlist}[(a)]
\item[(a)] $\Phi$ has \textup{bounded increments}: there exists $\xi
>0$ such that for {all $\tau$, we have} 
%
\[
\bigl|\Phi \bigl(\mathbf{X}(\tau+1) \bigr)-\Phi \bigl({\mathbf{X}(\tau)} \bigr)\bigr|\leq
\xi\qquad \mbox{almost surely};
\]
\item[(b)] \textup{Negative drift}: there exist $B>0$ and $\gamma>0$
such that
whenever\break  $\Phi(\mathbf{X}(\tau)) > B$,
\[
\mathbb{E} \bigl[\Phi \bigl(\mathbf{X}(\tau+1) \bigr)-\Phi \bigl(\mathbf{X}(\tau
) \bigr) | \mathbf{X}(\tau) \bigr]\leq-\gamma.
\]
%
\end{longlist}
Then, a stationary probability distribution $\pi$ exists,
and we have an exponential upper bound on the tail probability of $\Phi$
{under 
$\pi$:} for any $\ell\in\Zp$,
%
\begin{equation}
\label{eq:expbd} \mathbb{P}_{\pi} \bigl(\Phi(\mathbf{X})> B + 2\xi\ell
\bigr)\leq \biggl(\frac
{\xi}{\xi+\gamma} \biggr)^{\ell+1}.
\end{equation}
In particular, {in steady state,} all moments of $\Phi$ 
are finite, that is, for every $k \in\mathbb{N}$,
\[
\mathbb{E}_{\pi} \bigl[\Phi^k(\mathbf{X}) \bigr]<\infty.
\]
\end{theorem}

\newcommand{\bX}{\mathbf{X}}
Theorem \ref{thm:bgt01} is identical to
Theorem 1 in \cite{BGT01}, except that \cite{BGT01} imposed
the additional condition
$\mathbb{E}_{\pi}[\Phi(\mathbf{X})]<\infty$. However, the latter
condition is redundant. Indeed, using Foster--Lyapunov criteria (see
\cite{foss}, e.g.),
conditions (a) and (b) in Theorem \ref{thm:bgt01} imply that the
Markov chain $\bX$ has a unique
stationary distribution $\pi$. Furthermore, Theorem 2.3 in \cite
{Hajek-delay} establishes that under conditions (a)
and (b),
all moments of $\Phi(\mathbf{X})$ are finite in steady state.
We note that Theorem 2.3 in \cite{Hajek-delay} and Theorem 1 of \cite
{BGT01} provide the same qualitative information [exponential tail
bounds for
$\Phi(\mathbf{X})$]. However, \cite{BGT01} contains the more precise bound
(\ref{eq:expbd}), which we will use to prove Theorem
\ref{thm:inlimit} in Section \ref{sec:inlimit}.

\begin{pf*}{Proof of Theorem \ref{thm:sw2}}
{The finiteness of the moments follows immediately from the bound in
(\ref{eq:expbd}),
so we only prove the exponential bound (\ref{eq:expbd}).}
We apply Theorem \ref{thm:bgt01} to the Lyapunov function $L_{\alpha
}$ and the uniformized chain $\tilde{\bN}(\cdot)$.
Again, denote the stationary distribution of $\tilde{\bN}(\cdot)$ by
$\bpi$,
and note that this is also the unique stationary distribution of $\bN
(\cdot)$.
The proof consists of verifying conditions~(a) and~(b).

(a) \textit{Bounded increments}. We wish to show that with
probability $1$, there exists $\xi$ such that
\[
\bigl|L_{\alpha} \bigl(\tilde{\bN}(\tau+1) \bigr)-L_{\alpha} \bigl(
\tilde{ \bN}(\tau ) \bigr)\bigr|\leq\xi.
\]
As usual, write $\mathbf{n}= \tilde{\bN}(\tau)$ and $\mathbf
{n}+\bbm= \tilde{\bN}(\tau+1)$, then
$\bbm= \mathbf{0}$ or $\bbm= \pm\be_i$ for some $i\in\cI$
with probability $1$.
For $\alpha\geq1$,
\[
L_{\alpha}(\mathbf{n}) = \biggl[\sum_{i\in\cI}
w_i n_i^{\alpha
+1} \biggr]^{{1}/{(\alpha+1)}},
\]
and for $\alpha\in(0,1)$, by Lemma \ref{lm:prpt}(iii$'$), we have
\[
\sum_{i\in\cI}w_i n_i^{\alpha+1}
- 2\sum_{i\in\cI} w_i \leq (\alpha+1)\sum
_{i\in\cI}w_i H_{\alpha}(n_i)
\leq\sum_{i\in\cI}w_i n_i^{\alpha+1}
+ 2\sum_{i\in\cI} w_i.
\]
In general, for $r, s\geq0$ and $\beta\in[0,1]$,
%
\begin{equation}
\label{eq:normeq} (r+s)^{\beta}  \leq r^{\beta}+s^{\beta}.
\end{equation}
Thus, by inequality (\ref{eq:normeq}),
\begin{eqnarray*}
&&\biggl[\sum_{i\in\cI} w_i
n_i^{\alpha+1} \biggr]^{{1}/{(\alpha
+1)}} - \biggl[ 2\sum
_{i\in\cI} w_i \biggr]^{{1}/{(\alpha+1)}}\\
&&\qquad \leq
L_{\alpha}(\mathbf{n}) \leq \biggl[\sum_{i\in\cI}
w_i n_i^{\alpha+1} \biggr]^{{1}/{(\alpha
+1)}} + \biggl[
2\sum_{i\in\cI} w_i \biggr]^{{1}/{(\alpha+1)}}.
\end{eqnarray*}
Hence, for any $\alpha> 0$,
\begin{eqnarray*}
\bigl|L_{\alpha}(\mathbf{n}+\bbm)-L_{\alpha}(\mathbf{n})\bigr| &\leq& \biggl|
\biggl[\sum_{i\in\cI} w_i
(n_i+m_i)^{\alpha+1} \biggr]^{{1}/{(\alpha+1)}} -
\biggl[\sum_{i\in\cI} w_i
n_i^{\alpha
+1} \biggr]^{{1}/{(\alpha+1)}} \biggr|
\\
& &{} + 2 \biggl[ 2\sum_{i\in\cI} w_i
\biggr]^{{1}/{(\alpha+1)}}
\\
&\leq& \biggl[\sum_{i\in\cI} w_i
|m_i|^{\alpha+1} \biggr]^{{1}/{(\alpha+1)}} + 2 \biggl[ 2\sum
_{i\in\cI} w_i \biggr]^{{1}/{(\alpha+1)}}
\\
&\leq& \max_{i\in\cI}w_i^{{1}/{(\alpha+1)}} + 2 \biggl[ 2
\sum_{i\in\cI} w_i \biggr]^{{1}/{(\alpha+1)}},
\end{eqnarray*}
where the second last inequality follows from the triangle inequality.
Thus we can take $\xi= \max_{i\in\cI}w_i^{{1}/{(\alpha+1)}} +
2 [ 2\sum_{i\in\cI} w_i  ]^{{1}/{(\alpha+1)}}$,
which is a load-dependent constant.

(b) \textit{Negative drift}. The negative drift condition is
established in Theorem \ref{thm:sw-drift},
with $\gamma= \veps K$, for some positive load-dependent constant $K$.

Note that we have verified conditions (a) and (b) for the Lyapunov
function~$L_{\alpha}$.
To show the actual exponential probability tail bound for $\|\bN\|
_{\infty}$, note that
$L_{\alpha}(\bN)\geq K''\|\bN\|_{\infty}$, for some load-dependent
constant $K''$. By suitably
redefining the constants $B$, $\xi$ and $K$, the same form of
exponential probability tail bound
is established for $\|\bN\|_{\infty}$.
\end{pf*}

{
\section{\texorpdfstring{An important application: Interchange of limits ($\alpha=1$)}
{An important application: Interchange of limits (alpha = 1)}}\label{sec:inlimit}
In this section, we assume throughout that $\alpha= 1$ (the
proportionally-fair policy), and establish the validity of the heavy-traffic
approximation for networks in steady state.
We first provide the necessary preliminaries to
state our main theorem, Theorem \ref{thm:inlimit}.
In Section \ref{ssec:inlimit},
we state and prove Theorem \ref{thm:inlimit},
which is a consequence of Lemmas~\ref{lem:tight} and~\ref{lem:init_tight}.
Further definitions and background are provided in Section \ref
{ssec:technical},
along with the proofs of Lemmas \ref{lem:tight} and \ref{lem:init_tight}.
All definitions and background stated in this section are taken from
\cite{KW04} and \cite{KKLW09}.


\subsection{Preliminaries}\label{ssec:prelim}
We give a preview of the preliminaries that we will introduce
before stating Theorem \ref{thm:inlimit}.
The goal of this subsection is to
provide just enough background to be able to
state Theorem \ref{thm:diffapp},
the diffusion approximation result from \cite{KKLW09}.
To do this, we need a precise description of the process obtained in
the limit,
under the diffusion scaling.
This limiting process is a diffusion process,
called semimartingale reflecting Brownian motion (SRBM) (Definition
\ref{df:srbm}),
with support on a polyhedral cone. This polyhedral cone is defined
through the concept of an invariant manifold (Definition \ref{df:inv_mnfd}).

As in Section \ref{ssec:ssc}, we consider a sequence of networks
indexed by $r$,
where~$r$ is to be thought of as increasing to infinity along a
sequence. The incidence matrix~$\bA$,
the capacity vector $\bC$, and the weight vector $\bkappa$ do not
vary with $r$. Recall the \textit{heavy-traffic} condition---Assumption
\ref{asmp:fssc}, and the definitions of the workload $\mathbf{w}$,
the workload process $\bW$ and the lifting map $\Delta$ from
Definition \ref{df:lift}. We carry the notation from Section \ref{ssec:ssc},
so that $\btheta> \mathbf{0}$, and $\bnu^r \rightarrow\bnu>
\mathbf{0}$, $\bmu^r \rightarrow\bmu> \mathbf{0}$ and
$r(\bC- \bA\brho^r) \rightarrow\btheta$ as $r\rightarrow\infty$.
Recall that $\bA\brho= \bC$.
Let $\hat{\bN}^r$ and $\hat{\bW}^r$ be as
in (\ref{scaling}).

The continuity of the lifting map $\Delta$ will be useful in the sequel.
%
\begin{proposition}[(Proposition 4.1 in \cite{KKLW09})]\label{prop:lift-cont}
The function $\Delta\dvtx \Rpp^{|\cJ|} \rightarrow\Rpp^{|\cI|}$ is
continuous. Furthermore,
for each $\mathbf{w}\in\Rpp^{|\cJ|}$ and $c > 0$,
%
\begin{equation}
\label{eq:lift-inv} \Delta(c\mathbf{w}) = c\Delta(\mathbf{w}).
\end{equation}
\end{proposition}

\begin{definition}[(Invariant manifold)]\label{df:inv_mnfd}
A state $\mathbf{n}\in\Rpp^{|\cI|}$ is called \textit{invariant} if
$\mathbf{n}= \Delta(\mathbf{w})$,
where $\mathbf{w}= \bA\bM^{-1} \mathbf{n}$ is the workload,
and $\Delta$ the lifting map defined in Definition \ref{df:lift}.
The set of all invariant states is called the \textit{invariant manifold},
and we denote it by $\mathscr{M}$.
We also define the \textit{workload cone} $\mathscr{W}$ {by}
$\mathscr{W} = \bA\bM^{-1} \mathscr{M}$,
where $\bM= \diag(\bmu)$ is as defined in Definition \ref{df:lift}.
\end{definition}

{The} invariant manifold $\mathscr{M}$
is a polyhedral cone {and} admits an explicit characterization:
we can write it as
\[
\mathscr{M} = \biggl\{\mathbf{n}\in\Rpp^{|\cI|} \dvtx n_i =
\frac{\rho
_i (\mathbf{q}^T \bA)_i}{\kappa_i} \mbox{ for all } i \in\cI, \mbox{ for some } \mathbf{q}\in\Rpp
^{|\cJ|} \biggr\}.
\]

Denote the $j$th face of $\mathscr{M}$ by $\mathscr{M}^j$, which can
be written as
\begin{eqnarray*}
\mathscr{M}^j &\triangleq& \biggl\{\mathbf{n}\in\Rpp^{|\cI|}
\dvtx  n_i = \frac{\rho_i (\mathbf{q}^T \bA)_i}{\kappa_i} \mbox{ for all } i \in\cI,
\\
&&\hspace*{16pt}{} \mbox{for some } \mathbf{q}\in\Rpp^{|\cJ|} \mbox{ satisfying }
q_j = 0 \biggr\}.
\end{eqnarray*}
Similarly, denote the $j$th face of $\mathscr{W}$ by $\mathscr{W}^j$,
which can be written as
\[
\mathscr{W}^j \triangleq\bA\bM^{-1} \mathscr{M}^j.
\]

\subsubsection*{Semimartingale reflecting Brownian motion (SRBM)}
\begin{definition}\label{df:srbm}
Define the covariance matrix
\[
\bGamma= 2\bA\bM^{-1}\operatorname{diag}(\bnu)\bM^{-1}A^T.
\]
An \textit{SRBM} that lives in the cone $\mathscr{W}$,
has direction of reflection $\be_j$ {(the $j$th unit vector)} on the
boundary $\mathscr{W}^j$ for each $j \in\cJ$,
has drift $\btheta$ and covariance $\bGamma$ and has initial
distribution $\bseta_0$
on $\mathscr{W}$ is an adapted, $|\cJ|$-dimensional process $\hat
{\bW}(\cdot)$ defined on some
filtered probability space $(\Omega, \mathscr{F}, \{\mathscr{F}_t\},
\mathbb{P})$ such that:
\begin{enumerate}[(iii)]
\item[(i)] $\mathbb{P}$-a.s., $\hat{\bW}(t) = \hat{\bW}(0) + \hat
{\mathbf{X}}(t) + \hat{\mathbf{U}}(t)$ for all $t\geq0$;
\item[(ii)] $\mathbb{P}$-a.s., $\hat{\bW}(\cdot)$ has continuous
sample paths, $\hat{\bW}(t) \in\mathscr{W}$
for all $t\geq0$, and $\hat{\bW}(0)$ has initial distribution
$\bseta_0$;
\item[(iii)] under $\mathbb{P}$, $\hat{\mathbf{X}}(\cdot)$ is a
$|\cJ|$-dimensional Brownian motion starting at the origin with drift
$\btheta$
and covariance matrix $\bGamma$;
\item[(iv)] for each $j\in\cJ$, $\hat{U}_j(\cdot)$ is an adapted,
one-dimensional process such that \mbox{$\mathbb{P}$-a.s.},
\begin{enumerate}[(a)]
\item[(a)] $\hat{U}_j(0)=0$;
\item[(b)] $\hat{U}_j$ is continuous and nondecreasing;
\item[(c)] $\hat{U}_j(t) = \int_0^t \mathbb{I}_{\{\hat{\bW}(s)
\in\mathscr{W}^j\}}\,d\hat{U}_j(s)$
for all $t \geq0$.
\end{enumerate}
\end{enumerate}
The process $\hat{\bW}(\cdot)$ is called an SRBM with the data
$(\mathscr{W}, \btheta, \bGamma, \{\be_j \dvtx j\in\cJ\}, \bseta_0)$.
\end{definition}

\subsubsection*{Diffusion approximation for $\alpha=1$}

\begin{asmption}[(Local traffic)]\label{asp:local}
For each $j \in\cJ$, there exists at least one $i \in\cI$ such that
$A_{ji} > 0$ and
$A_{ki} = 0$ for all $k \neq j$.
\end{asmption}

{Under the local traffic condition, a diffusion approximation holds.}

\begin{theorem}[(Theorem 5.2 in \cite{KKLW09})]\label{thm:diffapp}
Assume that $\alpha= 1$ and that the local traffic condition,
Assumption \ref{asp:local}, holds.
Suppose that the limit distribution of $\hat{\bW}^r(0)$ as $r
\rightarrow\infty$ is
$\bseta_0$ (a probability measure on $\mathscr{W}$) and that
%
\begin{equation}
\label{cond:tight} \bigl\|\hat{\bN}^r(0)-\Delta \bigl(\hat{
\bW}^r(0) \bigr)\bigr\|_{\infty} \rightarrow0\qquad \mbox{in probability,
as } r \rightarrow\infty.
\end{equation}
Then, the distribution of $(\hat{\bW}^r(\cdot), \hat{\bN}^r(\cdot
))$ converges weakly
(on compact time intervals) %
as $r \rightarrow\infty$ to a continuous process $(\hat{\bW}(\cdot
), \hat{\bN}(\cdot))$, where
$\hat{\bW}(\cdot)$ is an SRBM with data $(\mathscr{W}, \btheta,
\bGamma, \{\be_j, j\in\cJ\}, \bseta_0)$
and $\hat{\bN}(t) = \Delta(\hat{\bW}(t))$ for all $t$.
\end{theorem}

\subsection{Interchange of limits} \label{ssec:inlimit}
We now know that for $\alpha= 1$, under the local traffic condition,
the diffusion approximation holds. That is, the scaled process $(\hat
{\bW}^r(\cdot), \hat{\bN}^r(\cdot))$
converges in distribution to $(\hat{\bW}(\cdot), \hat{\bN}(\cdot
))$, with $\hat{\bW}(\cdot)$ being an SRBM.
For any $r$, the scaled processes $\hat{\bN}^r(\cdot)$ also have
stationary distributions~$\bpi^r$,
since they are all positive recurrent. These results can be summarized
in the diagram that follows.
\[
\xymatrixcolsep{5pc} \xymatrix{ \hat{\bN}^r(\cdot) |_{[0,T]}
\ar[d]_{T \rightarrow\infty} \ar [r]^{r\rightarrow\infty}_{\mathrm{\normalfont Theorem\ \ref{thm:diffapp}}} & \hat{\bN}(\cdot)
|_{[0,T]}\ar[d]_{T \rightarrow\infty}^{?}
\\
\bpi^r \ar@{.>}[r]^{r\rightarrow\infty}_{?} & \hat{\bpi} }
\]

As can be seen from the diagram, two natural questions to ask are:
\begin{longlist}[(1)]
\item[(1)] Does the diffusion process $\hat{\bN}(\cdot)$ have a
stationary probability distribution, $\hat{\bpi}$?
\item[(2)] If $\hat{\bpi}$ exists and is unique, do the
distributions $\bpi^r$ converge to $\hat{\bpi}$?
\end{longlist}

Our contribution here is a positive answer to question (2). More specifically,
if $\hat{\bN}(\cdot)$ has a unique stationary probability
distribution $\hat{\bpi}$,
then $\bpi^r$ converges in distribution to $\hat{\bpi}$.

\begin{theorem}\label{thm:inlimit}
Suppose that $\alpha= 1$ {and} that the local traffic condition,
Assumption \ref{asp:local}, holds.
Suppose further that $\hat{\bN}(\cdot)$ has a unique stationary
probability distribution $\hat{\bpi}$.
For each $r$, let $\bpi^r$ be the unique stationary probability
distribution of $\hat{\bN}^r$. Then,
\[
\bpi^r \rightarrow\hat{\bpi}\qquad\mbox{in distribution, as } r
\rightarrow\infty.
\]
\end{theorem}

The line of proof of Theorem \ref{thm:inlimit} is fairly standard.
We first establish {tightness} of the set of distributions $\{\bpi^r\}
$ in Lemma \ref{lem:tight}.
Letting the processes $\hat{\bN}^r(\cdot)$ be initially distributed
as $\{\bpi^r\}$,
we translate this tightness condition into an initial condition
similar to (\ref{cond:tight}), in Lemma \ref{lem:init_tight}. We then
apply Theorem \ref{thm:diffapp}
to deduce the convergence of the processes $\hat{\bN}^r(\cdot)$,
which by stationarity, leads to the convergence of the distributions
$\bpi^r$.
We state Lemmas \ref{lem:tight} and \ref{lem:init_tight} below,
and defer their proofs to the next section.


\begin{lemma}\label{lem:tight}
Suppose that $\alpha= 1$. The set of probability distributions $\{\bpi
^r\}$ is tight.
\end{lemma}
%

\begin{lemma}\label{lem:init_tight}
Consider the stationary probability distributions $\bpi^r$ of $\hat
{\bN}^r(\cdot)$,
and let $\{\bpi^{r_k}\}$ be any convergent subsequence of $\{\bpi^r\}$.
Let $\hat{\bN}^r(0)$ be distributed as $\bpi^r$ for each $r$.
Then there exists a subsequence $r_{\ell}$ of $r_k$
such that\looseness=-1
%
\begin{equation}
\label{eq:init_tight} \bigl\|\hat{\bN}^{r_{\ell}}(0) - \Delta \bigl(
\hat{\bW}^{r_{\ell
}}(0) \bigr) \bigr\|_{\infty} \rightarrow0
\end{equation}\looseness=0
in probability as $\ell\rightarrow\infty$, that is, {such} that
condition (\ref{cond:tight})
holds for the subsequence $\{(\hat{\bW}^{r_{\ell}}(\cdot), \hat
{\bN}^{r_{\ell}}(\cdot))\}$.
\end{lemma}

\begin{pf*}{Proof of Theorem \protect\ref{thm:inlimit}}
Since $\{\bpi^r\}$ is tight by Lemma \ref{lem:tight},
Prohorov's theorem implies that $\{\bpi^r\}$ is relatively compact in
the weak topology.
Let $\{\bpi^{r_k}\}$ be a convergent subsequence of the set of
probability distributions $\{\bpi^r\}$,
and suppose that $\bpi^{r_k} \rightarrow\bpi$ as $k \rightarrow
\infty${,} in distribution.

Let $\hat{\bN}^r(0)$ be distributed as $\bpi^r$ for each $r$.
Then by Lemma \ref{lem:init_tight}, there exists a subsequence
$r_{\ell}$ of $r_k$
such that
\[
\bigl\|\hat{\bN}^{r_{\ell}}(0) - \Delta \bigl(\hat{\bW}^{r_{\ell
}}(0) \bigr)
\bigr\|_{\infty} \rightarrow0
\]
in probability as $\ell\rightarrow\infty$. Denote the distribution
of $\hat{\bW}^r(0)$ by $\bseta^r$.
Since $\bpi^{r_k} \rightarrow\bpi$ as $k \rightarrow\infty$,
$\bpi^{r_{\ell}} \rightarrow\bpi$ as $\ell\rightarrow\infty$ as well,
and $\bseta^{r_{\ell}} \rightarrow\bseta$ as $\ell\rightarrow
\infty$,
for some probability distribution $\bseta$.

We now wish to apply Theorem \ref{thm:diffapp} to the sequence $\{\hat
{\bN}^{r_{\ell}}(\cdot)\}$.
The only condition that needs to be verified is that $\bseta$ has
support on $\mathscr{W}$.
This can be argued as follows. Let $\hat{\bN}(0)$ have distribution
$\bpi$,
and let $\hat{\bW}(0) = \bA\bM^{-1} \hat{\bN}(0)$ be the
corresponding workload.
Then $\hat{\bW}^{r_{\ell}}(0) \rightarrow\hat{\bW}(0)$ in
distribution as $r \rightarrow\infty$,
and $\hat{\bW}(0)$ has distribution $\bseta$. The lifting map
$\Delta$ is continuous
by Proposition \ref{prop:lift-cont},
so $\Delta (\hat{\bW}^{r_{\ell}}(0) ) \rightarrow\Delta
(\hat{\bW}(0) )$
in distribution as $r \rightarrow\infty$. This convergence, together
with (\ref{eq:init_tight})
and the fact that $\hat{\bN}^{r_{\ell}}(0) \rightarrow\hat{\bN
}(0)$ in distribution,
{implies} that $\hat{\bN}(0)$ and $\Delta (\hat{\bW}(0)
)$ are identically distributed.
Now $\Delta (\hat{\bW}(0) )$ has support on $\mathscr{M}$,
so $\hat{\bN}(0)$ is supported on $\mathscr{M}$ as well,
and so $\hat{\bW}(0)$, hence $\bseta$, is supported on $\mathscr{W}$.

By Theorem \ref{thm:diffapp}, $(\hat{\bW}^{r_{\ell}}(\cdot), \hat
{\bN}^{r_{\ell}}(\cdot))$
converges in distribution to a continuous process $(\hat{\bW}(\cdot
), \hat{\bN}(\cdot))$.
Suppose that $\hat{\bW}(\cdot)$ and $\hat{\bN}(\cdot)$ have
unique stationary distributions $\hat{\bseta}$
and $\hat{\bpi}$, respectively. The processes $(\hat{\bW}^{r_{\ell
}}(\cdot), \hat{\bN}^{r_{\ell}}(\cdot))$ are stationary,
so $(\hat{\bW}(\cdot), \hat{\bN}(\cdot))$ is stationary as well.
Therefore, $\hat{\bW}(0)$ and $\hat{\bN}(0)$ are distributed as
$\hat{\bseta}$ and $\hat{\bpi}$, respectively.
Since $(\hat{\bW}^{r_{\ell}}(0), \hat{\bN}^{r_{\ell}}(0))
\rightarrow(\hat{\bW}(0), \hat{\bN}(0))$ in distribution,
we have that $\bseta^{r_{\ell}} \rightarrow\hat{\bseta}$ and
$\bpi^{r_{\ell}} \rightarrow\hat{\bpi}$ weakly as $\ell
\rightarrow\infty$.
This shows that $\bpi= \hat{\bpi}$ and $\bseta= \hat{\bseta}$.
Since $\{\bpi^{r_k}\}$ is an arbitrary convergent subsequence{,}
$\hat{\bpi}$ is the unique weak limit point of $\{\bpi^r\}$, and
this shows that
$\bpi^r \rightarrow\hat{\bpi}$ in distribution.
\end{pf*}

For Theorem \ref{thm:inlimit} to {apply}, we need to {verify} that
$\hat{\bN}(\cdot)$
[or equivalently, $\hat{\bW}(\cdot)$] has a unique stationary distribution.
The following theorem states that when $\kappa_i = 1$\vadjust{\goodbreak} for all $i\in
\cI$,
this condition holds; more specifically, the SRBM $\hat{\bW}(\cdot)$
has a unique stationary distribution, which turns out to have a product
form.

\begin{theorem}[(Theorem 5.3 in \cite{KKLW09})]\label{thm:srbm_dist}
Suppose that $\alpha= 1$ and $\kappa_i = 1$ for all $i \in\cI$. Let
$\hat{\bseta}$ be the measure
on $\mathscr{W}$ that is absolutely continuous with respect to
Lebesgue measure with density
given by
%
\begin{equation}
p(\mathbf{w}) = \exp \bigl(\langle\bv,\mathbf{w}\rangle \bigr),\qquad\mathbf {w}\in
\mathscr{W},
\end{equation}
where
%
\begin{equation}
\bv= 2\bGamma^{-1}\btheta.
\end{equation}
The product measure $\hat{\bseta}$ is an invariant measure for the
SRBM with state space~$\mathscr{W}$,
directions of reflection $\{\be_j, j\in\cJ\}$, drift $\btheta$ and
covariance matrix $\bGamma$.
After normalization, it defines the \textit{unique} stationary
distribution for the SRBM.
\end{theorem}

By {Theorems} \ref{thm:inlimit} and \ref{thm:srbm_dist},
the following corollary is immediate.
%
\begin{corollary}\label{cor:inlimit}
Suppose that $\alpha= 1$ and $\kappa_i = 1$ for all $i \in\cI$.
Suppose further that the local traffic condition, Assumption \ref
{asp:local}, holds.
Let $\hat{\bpi}$ be the unique stationary probability distribution of
$\hat{\bN}(\cdot)$.
For each $r$, let $\bpi^r$ be the unique stationary probability
distribution of $\hat{\bN}^r$. Then,
\[
\bpi^r \rightarrow\hat{\bpi}\qquad \mbox{in distribution, as } r
\rightarrow\infty.
\]
\end{corollary}

\subsection{\texorpdfstring{Proof of Lemmas \protect\ref{lem:tight} and \protect\ref{lem:init_tight}}
{Proof of Lemmas 7.7 and 7.8}}\label{ssec:technical}
\mbox{}
\begin{pf*}{Proof of Lemma \protect\ref{lem:tight}}
To establish tightness, it suffices to show that for every $y>0$ there
exists a
compact set $\mathbb{K}_y \subset\Rpp^{|\cI|}$
such that
%
\begin{equation}
\label{eq:tight1} \limsup_{r \rightarrow\infty} \bpi^r \bigl(
\Rpp^{|\cI|} \setminus \mathbb{K}_{{y}} \bigr)
\leq{e^{-y}}.
\end{equation}
We now proceed to define the compact {sets} $\mathbb{K}_{{y}}$.
As in the proof of Theorem \ref{thm:fssc}, let $\veps_r = \veps
(\brho^r)$ be the gap in the $r$th system.
Then, under Assumption \ref{asmp:fssc}, for sufficiently large $r$,
$\veps_r \geq D/r$ for some network-dependent constant $D>0$.
Since $\alpha= 1$, Theorem \ref{thm:sw2} implies that for the $r$th
system, there exist load-dependent constants
$K_r > 0$ and $\xi_r > 0$
such that for every $\ell\in\Zp$,
%
\begin{equation}
\label{eq:tight2} \mathbb{P}_{\bpi^r} \biggl(\bigl\|\bN^r
\bigr\|_{\infty} \geq\frac
{K_r}{\veps_r} + 2\xi_r \ell \biggr) \leq
\biggl(\frac{\xi_r}{\xi_r + \veps_r} \biggr)^{\ell+1}.
\end{equation}
{By the definition of a positive load-dependent constant,}
there exist continuous functions $f_1$ and $f_2$ on the open positive
orthant such that for all $r$,
$K_r = f_1(\bmu^r, \bnu^r)$ and $\xi_r = f_2(\bmu^r, \bnu^r)$.
Since $\bmu^r \rightarrow\bmu> \mathbf{0}$ and $\bnu^r \rightarrow
\bnu> \mathbf{0}$,
we have $K_r \rightarrow K \triangleq f_1(\bmu, \bnu) > 0$ and $\xi
_r \rightarrow\xi\triangleq f_2(\bmu, \bnu) > 0$.
Define
\[
\mathbb{K}_{{y}} \triangleq \biggl\{ \bv\in\Rpp^{|\cI|} \dvtx \|
\bv \| _{\infty} \leq\frac{(K+1) + 4(\xi+1)^2 \cdot{y}}{D} \biggr\}.\vadjust{\goodbreak}
\]
We now show that (\ref{eq:tight1}) holds, or equivalently, by the
definition of $\mathbb{K}_{{y}}$,
we show that for every $y > 0$,
%
\begin{equation}
\label{eq:tight3} \limsup_{r \rightarrow\infty} \PP_{\bpi^r} \biggl(
\frac{1}{r}\bigl\|\bN ^r\bigr\|_{\infty} > \frac{(K+1) + 4(\xi+1)^2 y}{D}
\biggr) \leq e^{-y}.
\end{equation}
Let $\ell_r \triangleq\lfloor2\xi_r y/\veps_r \rfloor$,
where for $z \in\R$, $\lfloor z \rfloor$ is the largest integer not
exceeding $z$.
By~(\ref{eq:tight2}), we have
\[
\PP_{\bpi^r} \biggl(\frac{1}{r}\bigl\|\bN^r
\bigr\|_{\infty} \geq\frac
{K_r}{r\veps_r} + \frac{2\xi_r\ell_r}{r} \biggr) \leq
\biggl(\frac{1}{1+{\veps_r}/{\xi_r}} \biggr)^{\ell_r+1}.
\]
Taking logarithms on both sides, we have
\[
\log\PP_{\bpi^r} \biggl(\frac{1}{r}\bigl\|\bN^r
\bigr\|_{\infty} \geq\frac
{K_r}{r\veps_r} + \frac{2\xi_r\ell_r}{r} \biggr) \leq-(
\ell_r + 1)\log \biggl(1+\frac{\veps_r}{\xi_r} \biggr).
\]
Since $\veps_r \rightarrow0$ and $\xi_r \rightarrow\xi> 0$ as $r
\rightarrow\infty$,
$\frac{\veps_r}{\xi_r} < 1$ for sufficiently large $r$.
Since $\log(1 + t) \geq t/2$ for $t \in[0, 1]$, we have
\[
-(\ell_r + 1)\log \biggl(1+\frac{\veps_r}{\xi_r} \biggr) \leq -(
\ell_r + 1) \frac{\veps_r}{2\xi_r},
\]
when $r$ is sufficiently large. By definition, $\ell_r = \lfloor2\xi
_r y/\veps_r \rfloor$,
so $\ell_r + 1 \geq2\xi_r x/\veps_r$, or equivalently,
$-(\ell_r + 1)\frac{\veps_r}{2\xi_r} \leq-y$.
Thus, when $r$ is sufficiently large,
\[
\log\PP_{\bpi^r} \biggl(\frac{1}{r}\bigl\|\bN^r
\bigr\|_{\infty} \geq\frac
{K_r}{r\veps_r} + \frac{2\xi_r\ell_r}{r} \biggr) \leq-y.
\]
Consider the term $\frac{K_r}{r\veps_r} + \frac{2\xi_r\ell_r}{r}$.
When $r$ is sufficiently large, $r \veps_r \geq D$, $K_r \leq K+1$,
and $\xi_r \leq\xi+1$,
and so
\[
\frac{K_r}{r\veps_r} + \frac{2\xi_r\ell_r}{r} \leq\frac
{K_r}{r\veps_r} +
\frac{2\xi_r (2\xi_r y)}{r\veps_r} \leq\frac{K+1}{D} + \frac{4(\xi+1)^2y}{D}.
\]
Thus, for sufficiently large $r$,
\begin{eqnarray*}
&& \log\PP_{\bpi^r} \biggl(\frac{1}{r}\bigl\|\bN^r
\bigr\|_{\infty} > \frac
{(K+1) + 4(\xi+1)^2 y}{D} \biggr)
\\
&&\qquad\leq \log\PP_{\bpi^r} \biggl(\frac{1}{r}\bigl\|\bN^r
\bigr\|_{\infty} \geq\frac{K_r}{r\veps_r} + \frac{2\xi_r\ell_r}{r} \biggr) \leq-y.
\end{eqnarray*}
This establishes (\ref{eq:tight3}), and also the tightness of $\{\bpi
^r\}$.
\end{pf*}

Next, we prove Lemma \ref{lem:init_tight}.
To this end, we need some definitions and background.
In particular, we need the concept and properties of \textit{fluid model
solutions}.

\begin{definition}
A \textit{fluid model solution} (FMS) is an absolutely continuous
function $\mathbf{n}\dvtx [0,\infty) \rightarrow\Rpp^{|\cI|}$
such that at each regular point\setcounter{footnote}{1}\footnote{A point $t \in(0,\infty)$
is a \textit{regular point}
of an absolutely continuous function $f \dvtx [0,\infty)\rightarrow\Rpp
^{|\cI|}$
if each component of $f$ is differentiable at $t$. Since
$\mathbf{n}$ is absolutely continuous, almost every time $t \in
(0,\infty)$ is a regular point for $\mathbf{n}$.}
$t>0$ of $\mathbf{n}(\cdot)$, we have, for each $i\in\cI$,
%
\begin{equation}
\label{eq:fluid} \frac{d}{dt}n_i(t) = \cases{ \nu_i - \mu_i \Lambda_i
\bigl(\mathbf{n}(t) \bigr), & \quad$\mbox{if } n_i(t)>0,$
\vspace*{2pt}\cr
0, & \quad$\mbox{if } n_i(t)=0,$}
\end{equation}
and for each $j\in\cJ$,
%
\begin{equation}
\label{eq:fluid2} \sum_{i\in\cI_{+}(\mathbf{n}(t))} A_{ji}
\Lambda_i \bigl(\mathbf{n}(t) \bigr) + \sum
_{i \in\cI_{0}(\mathbf{n}(t))} A_{ji}\rho_i \leq
C_j,
\end{equation}
where $\cI_{+}(\mathbf{n}(t)) = \{i\in\cI\dvtx  n_i(t)>0\}$ and $\cI
_{0}(\mathbf{n}(t)) = \{i\in\cI\dvtx  n_i(t)=0\}$.
Note that here $\bA\brho= \bC$.
\end{definition}

We now collect some properties of {a} FMS.
The following proposition states that the invariant manifold $\mathscr{M}$
consists exactly of the stationary {points of a} FMS.
%
\begin{proposition}[(Theorem 4.1 in \cite{KKLW09})]\label{prop:inv}
A vector $\mathbf{n}_0$ is an invariant state, that is, $\mathbf{n}_0
\in\mathscr{M}$,
if and only if for every fluid model solution $\mathbf{n}(\cdot)$
with $\mathbf{n}(0) = \mathbf{n}_0$, we have
$\mathbf{n}(t) = \mathbf{n}_0$ for all $t > 0$.
\end{proposition}

The following theorem states that starting from any initial condition,
a FMS will eventually be close to the invariant manifold $\mathscr{M}$.

\begin{theorem}[(Theorem 5.2 in \cite{KW04})]\label{thm:fluid1}
Fix $R\in(0,\infty)$ and $\delta> 0$. There is a constant
$T_{R,\delta} < \infty$ such
that for every fluid model solution $\mathbf{n}(\cdot)$ satisfying $\|
\mathbf{n}(0)\|_{\infty} \leq R$, we have
\[
d \bigl(\mathbf{n}(t), \mathscr{M} \bigr) < \delta \qquad\mbox{for all } t >
T_{R,\delta},
\]
where $d(\mathbf{n}(t), \mathscr{M}) \triangleq\inf_{\mathbf{n}\in
\mathscr{M}} \|\mathbf{n}- \mathbf{n}(t)\|_{\infty}$
is the distance from $\mathbf{n}(t)$ to the manifold~$\mathscr{M}$.
\end{theorem}

Proposition \ref{prop:lyap-dec} states that the value of the Lyapunov function
$F_1$ defined in (\ref{flpnv}) decreases along the path of any FMS.
%
\begin{proposition}[(Corollary 6.1 in \cite{KW04})]\label{prop:lyap-dec}
At any regular point $t$ of a fluid model solution $\mathbf{n}(\cdot
)$, we have
\[
\frac{d}{dt}F_1 \bigl(\mathbf{n}(t) \bigr) \leq0,
\]
and the inequality is strict if $\mathbf{n}(t) \notin\mathscr{M}$.
\end{proposition}

Using Proposition \ref{prop:lyap-dec} and the continuity of the
lifting map $\Delta$,
we can translate Theorem \ref{thm:fluid1} into the following version,
which will be used to prove Lemma \ref{lem:init_tight}.

\begin{lemma}\label{cor:fluid1}
Fix $R\in(0,\infty)$ and $\delta> 0$. There is a constant
$T_{R,\delta} < \infty$ such
that for every fluid model solution $\mathbf{n}(\cdot)$ satisfying $\|
\mathbf{n}(0)\|_{\infty} \leq R$ we have
\[
\bigl\|\mathbf{n}(t) - \Delta \bigl(\mathbf{w}(t) \bigr)\bigr\|_{\infty} < \delta\qquad
\mbox{for all } t > T_{R,\delta},
\]
where $\mathbf{w}(t) = \mathbf{w}(\mathbf{n}(t))$ is the workload
corresponding to $\mathbf{n}(t)$; see Definition~\ref{df:lift}.
\end{lemma}

\begin{pf}
Fix $R > 0$ and $\delta> 0$.
Let $\|\mathbf{n}(0)\|_{\infty} \leq R$. Then
\[
F_1 \bigl(\mathbf{n}(0) \bigr)= \frac{1}{2}\sum
_{i \in I} \nu_i^{-1}\kappa_i
n_i^2(0) \leq R',
\]
where $R'$ depends on $R$ and the system parameters.
Since $\mathbf{n}(\cdot)$ is absolutely continuous,
by Proposition \ref{prop:lyap-dec} and the fundamental theorem of calculus,
we have that
$F_1(\mathbf{n}(t))\leq R'$ for all $t \geq0$.
Define the set
\[
\cS\triangleq \bigl\{\mathbf{n}\in\Rpp^{|\cI|} \dvtx F_1(
\mathbf{n}) \leq R' \bigr\},
\]
and its $\delta$-\textit{fattening}
\[
\cS_{\delta} \triangleq \bigl\{\mathbf{n}\in\Rpp^{|\cI|} \dvtx \bigl\|
\mathbf {n}- \mathbf{n}'\bigr\| \leq\delta\mbox{ for some }
\mathbf{n}' \in\cS \bigr\}.
\]
Note that both $\mathcal{S}$ and $\mathcal{S}_{\delta}$ are compact
sets, and $\mathbf{n}(t)\in\mathcal{S} \subset\cS_{\delta}$ for
all $t \geq0$.

Now consider the workload $\mathbf{w}$ defined in Definition \ref{df:lift}.
Define the set $\mathbf{w}(\cS_{\delta}) = \{\bv\in\Rpp^{|\cJ|} \dvtx \bv= \mathbf{w}(\mathbf{n}) \mbox{ for some } \mathbf{n}\in\cS
_{\delta}\}$.
Since $\mathbf{w}$ is a linear map, there exists a load-dependent
constant $H$
such that
\[
\bigl\|\mathbf{w}(\mathbf{n})-\mathbf{w} \bigl(\mathbf{n}' \bigr)
\bigr\|_{\infty}\leq H\bigl\| \mathbf{n}-\mathbf{n}'\bigr\|_{\infty}
\]
for any $\mathbf{n}, \mathbf{n}' \in\Rpp^{|\cI|}$.
Thus $\mathbf{w}(\cS_{\delta})$ is also a compact set.
Since $\mathbf{n}(t) \in\cS_{\delta}$ for all $t \geq0$, $\mathbf
{w}(t) \in\mathbf{w}(\cS_{\delta})$
for all $t \geq0$.
By Proposition \ref{prop:lift-cont}, $\Delta$ is a continuous map, so
$\Delta$ is uniformly continuous when restricted to $\mathbf{w}(\cS
_{\delta})$.
Therefore, there exists $\delta' > 0$ such that for any $\mathbf{w}',
\mathbf{w}\in\mathbf{w}(\cS_{\delta})$ with
$\|\mathbf{w}' - \mathbf{w}\|_{\infty} < \delta'$, $\|\Delta
(\mathbf{w}') - \Delta(\mathbf{w})\|_{\infty} < \frac{\delta}{2}$.
Thus for any $\mathbf{n}, \mathbf{n}' \in\cS_{\delta}$ with $\|
\mathbf{n}- \mathbf{n}'\|_{\infty} < \delta'/H$, we have
$\|\mathbf{w}(\mathbf{n}) - \mathbf{w}(\mathbf{n}') \| \leq\delta
'$, and
\[
\bigl\|\Delta \bigl(\mathbf{w}(\mathbf{n}) \bigr) - \Delta \bigl(\mathbf{w} \bigl(
\mathbf {n}' \bigr) \bigr)\bigr\|_{\infty} < \frac{\delta}{2}.
\]
Let $\delta'' = \min\{\delta/2, \delta'/H\}$. By Theorem \ref{thm:fluid1},
there exists $T_{R,\delta''}$ such that for all $t \geq T_{R, \delta''}$,
\[
d \bigl(\mathscr{M},\mathbf{n}(t) \bigr) < \delta''.
\]
In particular, there exists $\mathbf{n}\in\mathscr{M}$ [which may
depend on $\mathbf{n}(t)$]
such that $\|\mathbf{n}- \mathbf{n}(t)\|_{\infty} < \delta'' <
\delta'/H$.
Since $\mathbf{n}(t) \in\cS$ and $\delta'' < \delta$, $\mathbf
{n}\in\cS_{\delta}$ as well. Thus
\[
\bigl\|\Delta \bigl(\mathbf{w}(\mathbf{n}) \bigr) - \Delta \bigl(\mathbf{w} \bigl(
\mathbf {n}(t) \bigr) \bigr)\bigr\|_{\infty} < \frac{\delta}{2}.
\]
By Proposition \ref{prop:inv}, since $\mathbf{n}\in\mathscr{M}$, we have
$\mathbf{n}= \Delta(\mathbf{w}(\mathbf{n}))$, and hence
\[
\bigl\|\mathbf{n}- \Delta \bigl(\mathbf{w} \bigl(\mathbf{n}(t) \bigr) \bigr)
\bigr\|_{\infty} < \frac
{\delta}{2}.
\]
Thus for all $t \geq T_{R, \delta''}$,
\begin{eqnarray*}
\bigl\|\mathbf{n}(t) - \Delta \bigl(\mathbf{w}(t) \bigr)\bigr\|_{\infty} &\leq& \bigl\|
\mathbf{n}- \mathbf{n}(t)\bigr\|_{\infty} + \bigl\|\mathbf{n}- \Delta \bigl(
\mathbf{w}(t) \bigr)\bigr\|_{\infty}
\\
& < & \delta'' + \frac{\delta}{2} \leq
\frac{\delta}{2}+\frac
{\delta}{2} = \delta.
\end{eqnarray*}
Note that $\delta''$ depends on $R$, $\delta$, and the system parameters.
Thus we can rewrite $T_{R, \delta''}$ as $T_{R, \delta}$. This
concludes the proof of the lemma.
\end{pf}

The last property of {a} FMS that we need is the {tightness of}
the fluid-scaled processes $\bar{\bN}^r$ and $\bar{\bW}^r$, defined by
%
\begin{equation}
\label{eq:fluidscale} \bar{\bN}^r(t) = \bN^r(rt)/r\quad
\mbox{and}\quad \bar{\bW}^r(t) = \bW^r(rt)/r.
\end{equation}

\begin{theorem}[(Theorem B.1 in \cite{KW04})]\label{thm:fluid2}
Suppose that $\{\bar{\bN}^r(0)\}$ converges in distribution as $r
\rightarrow\infty$
to a random variable taking values in $\mathbb{R}_+^{|\cI|}$.
Then the sequence $\{\bar{\bN}^r(\cdot)\}$ is $C$-tight,\footnote
{Consider the space $\mathbf{D}^{|\cI|}$ of functions $f\dvtx [0, \infty
) \rightarrow\R^{|\cI|}$
that are right-continuous on $[0, \infty)$ and have finite limits from
the left on $(0, \infty)$.
Let this space be endowed with the usual Skorohod topology; cf. Section
12 of \cite{billingsley1999convergence}.
The sequence $\{\bar{\bN}^r(\cdot)\}$ is \textit{tight} if the
probability measures induced on $\mathbf{D}^{|\cI|}$
are tight. The sequence is \textit{$C$-tight} if it is tight
and any weak limit point
{is a measure supported on the set of continuous sample paths.}}
and any weak limit point $\bar{\bN}(\cdot)$
of this sequence almost surely satisfies the fluid model equations
(\ref{eq:fluid})
and (\ref{eq:fluid2}).
\end{theorem}

\begin{pf*}{Proof of Lemma \ref{lem:init_tight}}
Consider the unique stationary distributions $\bpi^r$ of $\hat{\bN
}^r(\cdot)$,
and $\bseta^r$ of $\hat{\bW}^r(\cdot)$.
Let $\bpi^{r_k}$ be a convergent subsequence, and suppose that $\bpi
^{r_k} \rightarrow\bpi$
in distribution, as $k \rightarrow\infty$.
Suppose that at time $0$, $\frac{1}{r_k}\bN^{r_k}(0)$ is distributed
as $\bpi^{r_k}$.
Then $\frac{1}{r_k}\bW^{r_k}(0)$ is distributed as $\bseta^r_k$,
which converges in distribution as well, say to $\bseta$.

We now use the earlier stated FMS properties to prove the lemma.
Note that for all $r$,
\[
\frac{1}{r}\bN^r(0) = \bar{\bN}^r(0) = \hat{
\bN}^r(0) \quad\mbox{and}\quad \frac{1}{r}\bW^r(0) = \bar{
\bW}^r(0) = \hat {\bW}^r(0),
\]
and consider the fluid-scaled processes $\bar{\bN}^{r_k}(\cdot)$ and
$\bar{\bW}^{r_k}(\cdot)$.
Since $ \{\bar{\bN}^{r_k}(0) \}$ converges in distribution
to $\bpi$, Theorem \ref{thm:fluid2} implies that
the sequence $\{\bar{\bN}^{r_k}(\cdot)\}$ is $C$-tight, and any weak
limit $\bar{\bN}(\cdot)$ almost surely satisfies the fluid model equations.
Let $\bar{\bN}(\cdot)$ be a weak limit point of $ \{\bar{\bN
}^{r_k}(\cdot) \}$,
and suppose that the subsequence ${\{} \bar{\bN}^{r_{\ell}}(\cdot
){\}}$ of ${\{}\bar{\bN}^{r_k}(\cdot){\}}$ converges weakly to $\bar
{\bN}(\cdot)$.

Let $\delta> 0$. We will show that we can find $r(\delta)$ such that
for $r_{\ell} > r(\delta)$,
\[
\mathbb{P} \bigl(\bigl\|\bar{\bN}^{r_{\ell}}(0)-\Delta \bigl(\bar{\bW
}^{r_{\ell}}(0) \bigr)\bigr\|_{\infty} > \delta \bigr) < \delta.
\]
{Since $\bar{\bN}(0)$ is a well-defined random variable,}
there exists $R_{\delta} > 0$ such that
\[
\mathbb{P} \bigl(\bigl\|\bar{\bN}(0)\bigr\|_{\infty} > R_{\delta} \bigr) <
\frac{\delta}{2}.
\]
Now, for all sample paths $\omega$ such that $\|\bar{\bN}(0)(\omega
)\|_{\infty} \leq R_{\delta}$,
and such that $\bar{\bN}(\cdot)(\omega)$ satisfies the fluid model
equations,
Lemma \ref{cor:fluid1} implies that there exists $T \triangleq
T_{R_{\delta}, \delta}$ such that
\[
\bigl\|\bar{\bN}(T) (\omega) - \Delta \bigl(\bar{\bW}(T) \bigr) (\omega)
\bigr\|_{\infty
} < \delta.
\]
Since $\bar{\bN}(\cdot)$ satisfies the fluid model equations almost
surely, we have
\[
\mathbb{P} \bigl(\bigl\|\bar{\bN}(T) - \Delta \bigl(\bar{\bW}(T) \bigr)
\bigr\|_{\infty} < \delta \bigr) > 1 - \frac{\delta}{2}.
\]
Now for each $r$, $\bar{\bN}^r(0)$ is distributed according to the
stationary distribution $\bpi^r$,
so $\bar{\bN}^r(\cdot)$ is a stationary process.
Since $\bar{\bN}^{r_{\ell}}(\cdot) \rightarrow\bar{\bN}(\cdot)$
weakly as $\ell\rightarrow\infty$,
$\bar{\bN}$ is also a stationary process. 
Thus, $\bar{\bN}(T)$ and $\bar{\bN}(0)$ are both distributed
according to $\bpi$. This implies that
\[
\mathbb{P} \bigl(\bigl\|\bar{\bN}(0) - \Delta \bigl(\bar{\bW}(0) \bigr)
\bigr\|_{\infty} < \delta \bigr) > 1 - \frac{\delta}{2}.
\]
Furthermore, since $\bar{\bN}^{r_{\ell}}(0) \rightarrow\bar{\bN
}(0)$ in distribution,
\[
\mathbb{P} \bigl(\bigl\|\bar{\bN}^{r_{\ell}}(0) - \Delta \bigl(\bar{
\bW}^{r_{\ell
}}(0) \bigr)\bigr\|_{\infty} < \delta \bigr) \rightarrow
\mathbb{P} \bigl(\bigl\|\bar{\bN}(0) - \Delta \bigl(\bar{\bW}(0) \bigr)
\bigr\|_{\infty} < \delta \bigr)
\]
as $\ell\rightarrow\infty$.
Thus there exists $r(\delta)$ such that for all $r_{\ell} > r(\delta)$,
\[
\mathbb{P} \bigl(\bigl\|\bar{\bN}^{r_{\ell}}(0) - \Delta \bigl(\bar{
\bW}^{r_{\ell
}}(0) \bigr)\bigr\|_{\infty} < \delta \bigr) > 1 - \delta.
\]
Since $\delta> 0$ is arbitrary,
\[
\bigl\|\hat{\bN}^{r_{\ell}}(0) - \Delta \bigl(\hat{\bW}^{r_{\ell}}(0) \bigr)
\bigr\| _{\infty} = \bigl\|\bar{\bN}^{r_{\ell}}(0) - \Delta \bigl(\bar{
\bW}^{r_{\ell}}(0) \bigr)\bigr\| _{\infty} \rightarrow0,
\]
in probability.
\end{pf*}


\section{Conclusion}\label{sec:concl}
The results in this paper can be viewed from two different
perspectives. On the one hand, they provide much new
information on the qualitative behavior (e.g., finiteness of
the expected number of flows, bounds on steady-state tail probabilities
and finite-horizon maximum excursion probabilities, etc.) of
the $\alpha$-fair policies for bandwidth-sharing network models.
{At an abstract level,}
our results highlight
{the importance of relying on}
a suitable Lyapunov function. Even if a
network is shown to be stable by using a particular Lyapunov
function, different choices and more detailed analysis may lead to more powerful
bounds.
{At a more concrete level, we presented}
a generic method for deriving
full state space collapse from multiplicative state space collapse,
and {another} for deriving steady-state exponential tail bounds.
The methods and results in this paper can be extended to
general switched network models. Parallel results for a packet level
model are
detailed in \cite{STZ10}.

\section*{Acknowledgments}
We are grateful to Kuang Xu for reading the
manuscript and his comments.
We would also like to thank the reviewers and editors for constructive feedback,
which has greatly improved the readability of the manuscript.


%



\printaddresses

\end{document}